\documentclass[12pt]{article}

\usepackage{amsmath}
\usepackage{amssymb}
\usepackage{amsfonts}
\usepackage{a4}
\usepackage{graphicx}
\usepackage[abs]{overpic}
\usepackage{caption}
\usepackage{wrapfig}
\usepackage{float}
\usepackage{graphicx}

\newtheorem{theorem}{Theorem}[section]
\newtheorem{lemma}[theorem]{Lemma}

\newtheorem{corollary}[theorem]{Corollary}
\newtheorem{definition}[theorem]{Definition}

\newtheorem{rem}[theorem]{Remark}
\newtheorem{example}[theorem]{Example}
\newcommand{\Proof}{\par\noindent{\em Proof. }}
\newcommand{\eop}{\nopagebreak\hspace*{\fill}$\Box$\smallskip}

\newcommand{\N}{\Bbb N}

\newcommand{\R}{\Bbb R}

\def\Id{\mathbf{Id}}
\def\id{\mathbf{id}}

\def\eps{\varepsilon}

\def\e{\mathbf{e}}

\def\dist{\operatorname{dist}}
\def\Xint#1{\mathchoice
   {\XXint\displaystyle\textstyle{#1}}%
   {\XXint\textstyle\scriptstyle{#1}}%
   {\XXint\scriptstyle\scriptscriptstyle{#1}}%
   {\XXint\scriptscriptstyle\scriptscriptstyle{#1}}%
   \!\int}
\def\XXint#1#2#3{{\setbox0=\hbox{$#1{#2#3}{\int}$}
     \vcenter{\hbox{$#2#3$}}\kern-.5\wd0}}

\def\dashint{\Xint-}
\usepackage{color}

\numberwithin{equation}{section}

\newcommand{\BBB}{\color{black}} 
\newcommand{\EEE}{\color{black}}

\begin{document}

\begin{center}
\begin{Large}
{\bf {A derivation of linearized Griffith energies from nonlinear models}}
\end{Large}
\end{center}

\begin{center}
\begin{large}
Manuel Friedrich\footnote{Faculty of Mathematics, University of Vienna, 
Oskar-Morgenstern-Platz 1, 1090 Vienna, Austria. {\tt manuel.friedrich@univie.ac.at}}
\end{large}
\end{center}

\begin{center}
\today
\end{center}
\bigskip

\begin{abstract}
We derive Griffith functionals in the framework of linearized elasticity from nonlinear and frame indifferent energies in brittle fracture via $\Gamma$-convergence. The convergence is given in terms of rescaled displacement fields measuring the distance of deformations from piecewise rigid motions. The configurations of the limiting model consist of partitions of the material, corresponding piecewise rigid deformations and displacement fields which are defined separately on each component of the cracked body. Apart from the linearized Griffith energy the limiting functional comprises also the segmentation energy which is necessary to disconnect the parts of the specimen.
\end{abstract}
\bigskip

\begin{small}

\end{small}

\tableofcontents

%--------------------------------------------------------------------------
%--------------------------------------------------------------------------
\section{Introduction}
%--------------------------------------------------------------------------

A thorough understanding of crack formation in brittle materials is of great interest in both experimental sciences and theoretical studies. Starting with the seminal contribution by Francfort and Marigo \cite{Francfort-Marigo:1998}, where the displacements and crack paths are determined from an energy minimization principle, various variational models in the framework of free discontinuity problems have appeared in the literature over the past years. These so-called Griffith functionals comprising elastic and surface contributions generalize the original Griffith theory (see \cite{Griffith:1921}) which is based on the fundamental idea that the formation of fracture may be regarded as the competition of elastic bulk and surface energies.

For the sake of a simplified mathematical description the investigation of fracture models in the realm of linearized elasticity is widely adopted (see e.g.  \cite{Ambrosio-Coscia-Dal Maso:1997, Bellettini-Coscia-DalMaso:98, Bourdin-Francfort-Marigo:2008, Chambolle:2003, Chambolle:2004, Iurlano:13}) and has led to a lot of realistic applications in engineering as well as to efficient numerical approximation schemes (we refer to \cite{Bourdin:07, Bourdin-Francfort-Marigo:2000,  Ortner:13, Focardi-Iurlano:13, Negri:2003, Negri:06, SchmidtFraternaliOrtiz:2009} making no claim to be exhaustive). On the contrary, their nonlinear counterparts  are  usually significantly more difficult  to treat since in the regime of finite elasticity the energy density of the elastic contributions is genuinely geometrically nonlinear due to frame indifference rendering the problem highly non-convex. Consequently, in contrast to linear models already the fundamental question if minimizing configurations for given boundary data exist at all is a major difficulty. Even more challenging tasks in this context are the determination of the material behavior under expansion or compression, in particular the  derivation of specific cleavage laws. 

Consequently, for a deeper understanding of nonlinear models the identification of an effective linearized theory is desirable as in this way one may rigorously show  that in the small displacement regime the neglection of effects arising from the non-linearities is a good approximation of the problem. Moreover,  such a derivation is also interesting in the context of discrete systems. Previous investigations which were motivated by the analysis of cleavage laws for brittle crystals (see \cite{FriedrichSchmidt:2011, FriedrichSchmidt:2014.1} or the seminal paper \cite{Braides-Lew-Ortiz:06}) have shown that the most interesting regime for the elastic strains is given by $\sqrt{\eps}$, where $\eps$ denotes the typical interatomic distance. Consequently, a passage from discrete-to-continuum systems naturally involves a simultaneous linearization process. 

In elasticity theory the nonlinear-to-linear limit is by now well understood in various different settings via $\Gamma$-convergence (cf. \cite{Braides-Solci-Vitali:07, DalMasoNegriPercivale:02, Schmidt:08, Schmidt:2009}), where the passage is performed in terms of suitably rescaled  displacement fields measuring the distance of the deformation from a rigid motion and being the fundamental quantity on which the linearized elastic energy depends. In fracture mechanics, however, the relation between the deformation of a material and corresponding displacements is more complicated since the body may be disconnected by the jump set into various components. In fact, it turns out that, without passing to rescaled configurations, in the small strain limit  nonlinear Griffith energies converge to a limiting functional  which is finite for piecewise rigid motions and measures the \textit{segmentation energy} which is necessary to disconnect the body.

Obviously a major drawback of this simple limiting model appears to be the fact that it does not capture the elastic deformations which are typically present in the nonlinear models. Consequently, in order to arrive at a limiting model showing coexistence of elastic and surface contributions it is indispensable to pass to rescaled configurations similarly as in \cite{DalMasoNegriPercivale:02}. The goal of this article is to identify such an effective linearized Griffith energy as the $\Gamma$-limit of nonlinear and frame indifferent models in the small strain regime. To the best of our knowledge such a result has not yet been derived in the general setting of free discontinuity problems introduced by Ambrosio and De Giorgi \cite{DeGiorgi-Ambrosio:1988}.

The farthest reaching result in this direction seems to be a recent contribution by Negri and Toader \cite{NegriToader:2013} where a nonlinear-to-linear analysis is performed in the context of quasistatic evolution for a restricted class of admissible cracks. In particular, in their model the different components of the jump set are supposed to have a least positive distance rendering the problem considerably easier from an analytical point of view. In particular, the specimen cannot be separated into different parts  effectively leading to a simple relation between the deformation and the rescaled displacement field. On the other hand, in \cite{FriedrichSchmidt:2014.2} we have performed a simultaneous discrete-to-continuum and nonlinear-to-linear analysis for general crack geometries, but under the simplifying assumption that all deformations lie close to the identity mapping. 

In the present context we establish a limiting linearized Griffith functional  in a planar setting without  any a priori assumptions on the deformation and the crack geometry. We identify an effective model which appears to be more general than the energies which are widely investigated in the literature. Whereas in elasticity theory, in the approaches \cite{FriedrichSchmidt:2014.2, NegriToader:2013}  mentioned before and in most linear fracture models there is a simple relation between the deformation of the material and the associated infinitesimal displacement field, in our framework  the deformation is related to a triple consisting of a partition of the domain, a corresponding piecewise rigid motion being constant on each connected component of the cracked body and a displacement field which is defined separately on each piece of the specimen. 

On each component of the partition the energy is of Griffith-type in the realm of linearized elasticity. In addition, the functional contains the segmentation energy  which is necessary to disconnected the parts of the body. In particular, the latter contribution  is a specific feature of our general model where we do not restrict the analysis to a linearization around a fixed rigid motion. 

 Let us briefly note that although all arguments used in the proofs of this article are valid in any space dimension, we have to restrict our analysis to two dimensions as one of the ingredients of our analysis, an $SBD$-rigidity result (see \cite{Friedrich-Schmidt:15}), has only been derived in a planar setting for isotropic surface energies. However, we believe that the estimate in \cite{Friedrich-Schmidt:15}  may be generalized in the future and then the generalizations for the results in the work at hand immediately follow.

As applications of our result \BBB we investigate problems with external forces and also \EEE present a cleavage law in a continuum setting with isotropic surface energies. As discussed before, the identification of critical loads and the investigation of crack paths is a challenging problem particularly for nonlinear models. The arguments in \cite{FriedrichSchmidt:2011, FriedrichSchmidt:2014.1, Mora:2010}, where boundary value problems of uniaxial extension for brittle materials were investigated, fundamentally relied on the application of certain slicing techniques and due to the lack of convexity were not adapted to treat the case of  compression. Our general $\Gamma$-limit result can now be applied to solve also boundary value problems of uniaxial compression which is as  the  uniaxial tension test a natural and interesting problem. Hereby we may complete the  picture about the derivation of cleavage laws in \cite{FriedrichSchmidt:2011, FriedrichSchmidt:2014.1}.

One essential point in our investigation is the establishing of a compactness result providing limiting configurations which consist of piecewise rigid motions and corresponding displacement fields. Similarly as in the derivation of linearized systems for elastic materials (see e.g. \cite{DalMasoNegriPercivale:02}), where the  main ingredient is a quantitative geometric rigidity estimate by Friesecke, James and M\"uller \cite{FrieseckeJamesMueller:02}, the starting point of our analysis is a quantitative $SBD$-rigidity result (see \cite{Friedrich-Schmidt:15}) in the framework of special functions of bounded deformation (see \cite{Ambrosio-Coscia-Dal Maso:1997, Bellettini-Coscia-DalMaso:98}),  which is  tailor-made for general Griffith models with coexistence of both both energy forms.

 As there is no uniform bound on the functions, it turns out that the limiting displacements are generically not summable and we naturally end up in the space of $GSBD$ functions (for the definition and basic properties we refer to \cite{DalMaso:13}). We believe that our results are interesting also outside of this specific context as they allow to solve more general variational problems in fracture mechanics. Typically, for compactness results in function spaces as $SBV$ (see \cite{Ambrosio-Fusco-Pallara:2000} for the definition and basic properties) and $SBD$ one needs $L^\infty$ or $L^1$ bounds on the functions (see \cite{Ambrosio:90, Bellettini-Coscia-DalMaso:98, DalMaso:13}). However, in many applications, in particular for atomistic systems and for models dealing with rescaled deformations, such estimates cannot be inferred from energy bounds. Nevertheless, we are able to treat problems without any a priori bound by passing from the deformations to displacement fields whose distance from  rigid motions can be controlled. 

The other essential point in our analysis is the investigation of the limiting configurations. In particular, we study the properties of the partition which disconnects the body into various parts. It turns out that an even finer segmentation may occur if on a connected component of the partition  the jump set of the corresponding displacement field further separates the body. Here it becomes apparent that we treat a real multiscale model as the jump heights at the boundaries associated to the coarse partition are of order $\gg  \sqrt{\eps}$ ($\sqrt{\eps}$ denotes the regime of the typical elastic strain), whereas the jump heights of the finer partition are of order $\sqrt{\eps}$. Moreover, it is evident that the choice of the limiting partition is not unique. However, we propose a selection principle and show existence and uniqueness of a \emph{coaresest partition}.

The paper is organized as follows. In Section \ref{rig-sec: main} we state the main compactness and $\Gamma$-convergence results and discuss properties of the limiting linearized  Griffith functional. Moreover, we present our application to cleavage laws for uniaxially extended or compressed brittle materials. 

Section \ref{rig-sec: pre} is devoted to some preliminaries. We first give the definition of special functions of bounded variation and deformation and discuss basic properties. Afterwards, we recall the notion of Caccioppoli partitions which will be fundamental in our analysis to analyze the properties of limiting configurations. Moreover, we recall geometric rigidity results for elastic and brittle materials, in particular the $SBD$-rigidity result proved in \cite{Friedrich-Schmidt:15}. 
  
In Section \ref{rig-sec: sub, comp1} we then establish the main compactness result for a sequence of deformations $(y_\eps)_\eps$, where $\eps$ stands for the order of the elastic energy. First, the convergence of the partitions and the corresponding rigid motions is based on compactness theorems for Caccioppoli partitions and piecewise constant functions (see \cite{Ambrosio-Fusco-Pallara:2000} or Section \ref{rig-sec: sub, cacciop} below).

Although the $SBD$-rigidity estimate is a fundamental ingredient in our analysis giving $L^2$ bounds for rescaled displacement fields, we still have to face major difficulties since the rigidity estimate provides a family of displacement fields $(u^\rho_\eps)^\rho_\eps$ with an additional parameter $\rho$ representing a `modification error' between $y_\eps$ and $u^\rho_\eps$. Consequently, the goal will be to choose an appropriate diagonal sequence.

An additional challenge is the fact that the bounds in the $SBD$-rigidity estimate depend on $\rho$ and blow up for $\rho \to 0$. For the symmetric part of the gradient this problem can be bypassed by a Taylor expansion taking the nonlinear elastic energy $\eps$ and a higher order term into account, which shows that the constant may be chosen independently of $\rho$. For the function itself, however, the problem is more subtle since a uniform bound cannot be inferred by energies bounds. In particular, generically the limiting configurations are not in $L^2$, but only  finite almost everywhere. The strategy to establish the latter assertion is to show that for fixed $\eps$ the functions $(u^\rho_{\eps})_\rho$ essentially coincide in a certain sense on the bulk part of the domain. Afterwards, by a careful analysis we can derive that such a property is preserved in the limit $\eps \to 0$, whereby we can establish a kind of equi-integrability of the configurations.

In Section \ref{rig-sec: sub, comp2} we concern ourselves with the limiting configurations consisting of a partition, a corresponding piecewise rigid motion and a displacement field. Recalling that genuinely the limits provided by the compactness result are highly non-unique we introduce the notion of a \textit{coarsest partition}. Roughly speaking, the definition states that the jump heights at the boundaries associated to this  partition are of order $\gg  \sqrt{\eps}$ leading to a meaningful mathematical description of the observation that the size of the crack opening is a multiscale phenomenon in our model. 

The fundamental point is the proof of existence and uniqueness of the coarsest partition. Uniqueness follows from the fact  that under the assumption that there are two different coarsest partitions one always can find an even coarser partition. Existence is a more challenging problem. We first give an alternative characterization and identify coarsest partitions as the maximal elements of the partial order on the set of admissible
partitions which is induced by subordination. We then show that each chain of the partial order has an upper bound repeating some arguments of the main compactness result. Consequently, the claim is inferred by an application of Zorn's lemma. Finally, having found the coarsest partition we can then show that the corresponding admissible displacement field is uniquely determined up to piecewise infinitesimal rigid motions. 

In Section \ref{rig-sec: sub, gamma1} we derive the main $\Gamma$-limit, where the elastic part can be treated as in \cite
{FrieseckeJamesMueller:02, Schmidt:2009} and for the surface energy we separate the effects arising from the segmentation energy and the crack energy inside the components by employing a structure theorem for Caccioppoli partitions (see Theorem \ref{th: local structure} below). \BBB At this point we also establish a result including external loads. \EEE

Finally, in Section \ref{rig-sec: sub, gamma2} we prove a cleavage law and extend the results obtained in \cite{FriedrichSchmidt:2011, FriedrichSchmidt:2014.1, Mora:2010} to the case of uniaxial compression, where we essentially follow the proof in \cite{FriedrichSchmidt:2014.2, Mora:2010}, in particular using a piecewise rigidity result in  $SBD$  (see \cite{Chambolle-Giacomini-Ponsiglione:2007}) and a structure theorem on the boundary of sets of finite perimeter (see \cite{Federer:1969}). It turns out that in the linearized limit the behavior for compression and extension is virtually identical. We briefly note that to avoid unphysical effects such as self-penetrability further modeling assumptions would be necessary.

%----------------------------------------------------------------------------------------------------------------------
\section{The model and main results}\label{rig-sec: main}
%-----------------------------------------------------------------------------------------------------------------------

\subsection{The nonlinear model}

Let $\Omega \subset \R^2$ open,  bounded with Lipschitz boundary. Recall the properties of the space $SBV(\Omega,\R^2)$, frequently abbreviated as $SBV(\Omega)$ hereafter, in Section \ref{rig-sec: sub, sbv}. Fix a (large) constant $M>0$ and define 
\begin{align}\label{rig-eq: SBVfirstdef}
SBV_M(\Omega) = \Big\{ y \in SBV(\Omega,\R^2):  \Vert y\Vert_{\infty} + \Vert \nabla y\Vert_{\infty} \le M, \ {\cal H}^1(J_y) < + \infty \Big\}.
\end{align}
Let $W:\R^{2 \times 2} \to [0,\infty)$ be a frame-indifferent stored energy density with $W(F) = 0$ iff $F \in SO(2)$. Assume that $W$ is continuous, $C^3$ in a neighborhood of $SO(2)$ and scales quadratically at $SO(2)$ in the direction perpendicular to infinitesimal rotations. In other words, there is a positive constant $c$ such that
\begin{align}\label{eq:W}
W(F) \ge c\dist^2(F,SO(2)) \ \ \  \BBB \text{for all} \ \ \ F \in \R^{2 \times 2} \ \ \ \text{with} \ \ \ |F|\le M. \EEE
\end{align}For $\eps >0$ define the Griffith-energy $E_\eps : SBV_M(\Omega) \to [0,\infty)$ by
\begin{align}\label{rig-eq: Griffith en}
E_\eps(y) =  \frac{1}{\eps}\int_\Omega W(\nabla y(x)) \,dx + {\cal H}^1(J_y).
\end{align}
We briefly note that we can also treat inhomogeneous materials where the energy density has the form $W: \Omega \times \R^{2 \times 2} \to [0,\infty)$. Moreover, it suffices to assume $W \in C^{2,\alpha}$, where $C^{2,\alpha}$ is the H\"older space with exponent $\alpha >0$.

The main goal of the present work is the identification of an effective linearized Griffith energy in the small strain limit which is related to the nonlinear energies $E_\eps$ through $\Gamma$-convergence. \BBB In this context, we also discuss minimization problems associated to $E_\eps$ for given body forces or boundary data. \EEE Moreover, we will investigate the limiting model which appears to be more general than many other Griffith functionals in the realm of linearized elasticity (cf. e.g. \cite{Bourdin-Francfort-Marigo:2008, Chambolle:2003, Chambolle:2004, Iurlano:13, SchmidtFraternaliOrtiz:2009}) as the limiting configuration not only consists of a displacement field, but also of a coarse partition of the domain and associated rigid motions. In particular, it will turn out that there are various scales for the size of the crack opening occurring in the system.

\BBB

\begin{rem}\label{rem:M}

{\normalfont The threshold $M$ in \eqref{rig-eq: SBVfirstdef} may be chosen arbitrarily large, but is fixed. Confining $y$ in this way effectively models a large box containing
the deformed specimen. The restriction on $\nabla y$ is necessary for technical reasons as it allows us to
apply a quantitative piecewise rigidity estimate, see Theorem \ref{rig-th: rigidity}.

Let us mention that (almost) minimizers of the nonlinear energy $\int_{\Omega\setminus \overline{J_y}}W(\nabla y)$ (for given boundary data) are possibly not Lipschitz continuous as particularly at nonsmooth points of the boundary $\partial (\Omega\setminus \overline{J_y})$ (e.g. at crack tips) the deformation gradient is expected to form singularities. Consequently, the constraint $\Vert \nabla y\Vert_{\infty} \le M$ is a real restriction on the class of admissible configurations from a mathematical point of view.

On the other hand, for materials undergoing brittle fracture there is typically a critical strain (and stress), beyond which failure occurs, and therefore the uniform bound on the absolute continuous part of the gradient  has a reasonable mechanical interpretation. Moreover, the energy of certain atomistic systems can be related to \eqref{rig-eq: Griffith en}  when deformations are identified with piecewise affine interpolations on cells of microscopic size (see e.g. \cite{Braides-Gelli:2002-2, FriedrichSchmidt:2014.2}). (Note that in discrete systems the parameter $\eps$ represents not only the order of the elastic energy, but also the typical interatomic distance.) In this context, the bound $\Vert \nabla y\Vert_{\infty} \le M$ is naturally satisfied. In fact, on cells exceeding such a threshold, also called \emph{ultimate strain} (see \cite{Braides-DalMaso-Garroni:1999}), a  discontinuous interpolation with bounded deformation gradient is introduced and their contribution to the energy then enters through the surface part of the energy functional.

Finally, let us mention that particularly from a computational point of view it is interesting to combine a continuum model as \eqref{rig-eq: Griffith en} with an atomistic approach  using the quasicontinuum method introduced in \cite{Tadmor}. Here the underlying idea is to split the domain into a bulk part with  a coarse, continuum description, and into certain \emph{critical regions} characterized by fast variations of the deformation gradient (such as regions near a dislocation core or a crack tip) where the problem is treated as a fully atomistic system at scale $\eps$ (see \cite{Miller}).  
}
\end{rem}

\EEE

\subsection{The segmentation problem}\label{rig-sec: sub, seg}

As a first natural approach to  the problem we concern ourselves with the question if the functionals $E_\eps$ can be related to a limiting functional for $\eps \to 0$ in terms of the deformations. We observe that for configurations with uniformly bounded energy $E_\eps(y_\eps)$ the absolute continuous part of the gradient satisfies $\nabla y_\eps \approx SO(2)$ as the stored energy density is frame-indifferent and minimized on $SO(2)$. Assuming that $y_\eps \to y$ in $L^1$, one can show that $\nabla y \in SO(2)$ a.e. applying lower semicontinuity results for $SBV$ functions (see \cite{Kristensen:1999}) and the fact that the quasiconvex envelope of $W$ is minimized exactly on $SO(2)$ (see \cite{Zhang:2004}). 

A piecewise rigidity result by Chambolle, Giacomini and Ponsiglione (see Theorem \ref{rig-cor: cgp} below) generalizing the classical Liouville result for smooth functions now states that an $SBV$ function $y$ satisfying the constraint $\nabla y \in SO(2)$ a.e. is a collection of an at most countable family of rigid deformations, i.e.  the body may be divided into different components each of which subject to a different rigid motion. 

Consequently, the limit of the sequence $E_\eps$ (in the sense of $\Gamma$-convergence) is given by the functional which is finite for piecewise rigid motions and measures the \textit{segmentation energy} which is necessary to disconnect the body. The exact statement is formulated in Corollary \ref{rig-cor: gamma} as a direct consequence of our main  $\Gamma$-convergence result in Theorem \ref{rig-th: gammaconv}.

\BBB Apparently  this simple limiting model does not account for \EEE the elastic deformations which are typically present in the nonlinear models. Consequently, to obtain a better understanding of the problem  it is desirable to pass to rescaled configurations and to derive a limiting linearized energy as it was performed in \cite{DalMasoNegriPercivale:02} in the framework of nonlinear elasticity theory. The main ingredient in  that analysis is  a quantitative rigidity result due to Friesecke, James and M\"uller (see Theorem \ref{rig-th: geo rig}). The starting point for our analysis will be a corresponding quantitative result in the $SBD$ setting (see \cite{Friedrich-Schmidt:15} or Theorem \ref{rig-th: rigidity}) adapted for Griffith functionals of the form \eqref{rig-eq: Griffith en} where  both elastic bulk and surface contributions are present.

%----------------------------------------------------------------------------------------------------------------------
\subsection{Compactness \BBB and limiting configurations\EEE}\label{rig-sec: sub, main com}
%-----------------------------------------------------------------------------------------------------------------------

We now present our main compactness result for rescaled displacement fields. As a preparation, recall the notion and basic properties of a \textit{Caccioppoli partition} in Section \ref{rig-sec: sub, cacciop}. For a given (ordered) Caccioppoli partition ${\cal P} = (P_j)_j$ of $\Omega$ let 
\begin{align}\label{rig-eq: defA}
{\cal R}({\cal P}) = \Big\{ T: \Omega \to \R^2:  \ T(x) = \sum\nolimits_j \chi_{P_j}(x) (R_j \, x + b_j), \  R_j \in SO(2), \  b_j \in \R^2  \Big\}
\end{align}
be the set of corresponding piecewise rigid motions. Likewise we define the set of piecewise infinitesimal rigid motions, denoted by ${\cal A}({\cal P})$, replacing $R_j \in SO(2)$ by $A_j \in \R^{2 \times 2}_{\rm skew} = \lbrace A \in \R^{2 \times 2}: A=-A^T\rbrace$. Moreover, we define the triples
\begin{align}\label{eq:triples}
 {\cal D} &:= \big\{ (u,{\cal P}, T): \ u \in   SBV(\Omega),   \ {\cal P} \text{ C.-partition of } \Omega,\  T \in {\cal R}({\cal P})  \big\}, \\
 {\cal D}_\infty &:= \big\{ (u,{\cal P}, T): \, {\cal P} \text{ C.-partition of } \Omega,\,  T \in {\cal R}({\cal P}), \, \nabla T^T u \in GSBD^2(\Omega)  \big\}.\notag
\end{align}
Here $\nabla T$ denotes the absolutely continuous part of $DT$. The space $GSBD^2(\Omega)$  generalizes the definition of the space $SBD(\Omega)$ based on certain slicing properties, see Section \ref{rig-sec: sub, sbv}. Define $e(G) = \frac{G^T + G}{2}$ for all $G \in \R^{2 \times 2}$ and denote by $\partial^*$ the \emph{essential boundary} (see below \eqref{eq: essential boundary}). \BBB Let $A \triangle B$ be the symmetric difference of two sets $A,B \subset \R^2$. \EEE We now formulate the main compactness theorem.

\begin{theorem}\label{rig-th: comp1}
Let $\Omega \subset \R^2$ open, bounded with Lipschitz boundary. Let $M>0$ and $\eps_k \to 0$ as $k \to \infty$. If $E_{\eps_k}(y_{k}) \le C$ for a sequence $(y_k)_k \subset SBV_M(\Omega)$, then there exists a subsequence (not relabeled) such that the following holds: \\
There are triples $(u_k, \BBB {\cal P}_k, \EEE T_k) \in {\cal D}$, where ${\cal P}_k = (P^k_j)_j$, \BBB and $c>0$ with \EEE
\begin{align}\label{rig-eq: comp2}
\begin{split}
(i)& \ \ u_k(x) - \eps_k^{-1/2} (y_k(x)  -T_k(x)) \to 0 \ \text{for a.e. } x \in \Omega \ \ \text{ for $k \to \infty$},\\
(ii)& \ \  \Vert\nabla u_k\Vert_{L^\infty(\Omega)} \le  c\eps_k^{-1/8}  \ \ \text{ for $k \in \N$}
\end{split}
\end{align}
  such that we find a limiting triple $(u, {\cal P}, T) \in {\cal D}_\infty$ with   
\begin{align}\label{rig-eq: comp1}
\begin{split}
(i)& \ \  \BBB |P^k_j \triangle P_j| \to 0 \EEE  \ \ \ \text{ for all} \ j \in \N,\\
(ii)& \ \   T_k \to T \text{ in } L^2(\Omega, \R^2), \ \ \ \nabla T_k \to \nabla T \text{ in } L^2(\Omega, \R^{2 \times 2})
\end{split}
\end{align}
for $k \to \infty$ and 
\begin{align}\label{rig-eq: comp1-2}
\begin{split}
(i)&\ \ u_k \to u \ \ \text{ a.e. in } \  \Omega,\\
(ii)& \ \ e ( \nabla T^T_k \nabla   u_k   ) \rightharpoonup  e (\nabla T^T \nabla u)  \ \ \text{ weakly in} \ L^2(\Omega,\R^{2\times 2}_{\rm sym})
\end{split}
\end{align}
for $k \to \infty$. Moreover, \BBB for the elastic and crack \EEE  energy we obtain  
\begin{align}\label{rig-eq: comp3}
(i)&\ \   \frac{1}{\eps_k}\int_\Omega W(\nabla y_k) + o(1) \ge \frac{1}{\eps_k}\int_\Omega   W(\Id + \sqrt{\eps_k} \nabla T^T_k \nabla   u_k  ) \ \ \ \text{as} \ k\to\infty,\notag\\
(ii)&\ \ \liminf_{k \to \infty} {\cal H}^1(J_{y_k}) \ge  {\cal H}^1\Big(\bigcup\nolimits_j \partial^* P_j \cap \Omega\Big) + {\cal H}^1\Big(J_u \setminus \bigcup\nolimits_j \partial^* P_j\Big).
\end{align}
\end{theorem}

\BBB
Recall that the central object in linearized elasticity is the symmetric part of the gradient, which comes from the fact that (1) deformations are linearized around the identity and (2) the orthogonal space to $SO(d)$ at the identity is given by the symmetric matrices (see e.g.  \cite{DalMasoNegriPercivale:02}). In the present context, where we possibly linearize around different rigid motions, the symmetrized gradient is accordingly replaced by $e(\nabla T^T u)$  in both the limiting description and the convergence (see \eqref{eq:triples} and \eqref{rig-eq: comp1-2}(ii), respectively).

In \eqref{rig-eq: comp1} and \eqref{rig-eq: comp1-2} the convergence for the partitions, rigid motions and displacement fields is given, respectively. Moreover, \eqref{rig-eq: comp2} and \eqref{rig-eq: comp3} represent \emph{compatibility conditions} for the triple $(u_k,{\cal P}_k,T_k)$: In general, $u_k$ is a modification of the rescaled displacement $\eps_k^{-1/2} (y_k -T_k)$, but asymptotically both configurations coincide (see \eqref{rig-eq: comp2}(i)).  Moreover, the modifications can be constructed such that $\nabla u_k$ is suitably controllable. (The exponent $-\frac{1}{8}$ is chosen for definiteness only and could be replaced by any small negative exponent, cf.  \eqref{rig-eq: main properties2}(iv) and the paragraph below Theorem \ref{rig-th: rigidity}.) Finally, the elastic and crack energy associated to the triples are controlled by the corresponding energies of $y_k$ up to small errors vanishing in the limit (see \eqref{rig-eq: comp3}).

\begin{definition}\label{def:conv}
{\normalfont
We say a sequence $(y_k)_k \subset SBV_M(\Omega)$ is \emph{asymptotically represented} by a limiting triple $(u,{\cal P},T) \in {\cal D}_\infty$, and write $y_k \to (u,{\cal P},T)$, if there is a sequence of triples $(u_k,{\cal P}_k,T_k) \in {\cal D}$ such that \eqref{rig-eq: comp2}-\eqref{rig-eq: comp3} are satisfied.
}
\end{definition}
Although we use the notation $\to$ and call $(u,{\cal P},T)$ a limiting triple, it is clear that Definition \ref{def:conv} cannot be understood as a convergence in the usual sense. In particular, in the small strain limit a tripling of the variables occurs, which is a specific feature of our limiting model. Additionally, the  triples  $(u,{\cal P},T)$ given by the main compactness theorem for a sequence $(y_k)_k$ are not determined uniquely, but crucially depend on the choice of the sequences $({\cal P}_k)_k$ and $(T_k)_k$. To illustrate the latter phenomenon, \EEE we consider the following example.

\begin{example}\label{ex}
{\normalfont
Consider $\Omega= (0,3) \times (0,1)$, $\Omega_1 = (0,1) \times (0,1)$, $\Omega_2 = (1,2) \times (0,1)$, \BBB $\Omega_3 = (2,3) \times (0,1)$ \EEE and   
$$y_k = \id \chi_{\Omega_1} + (\id + \alpha\sqrt{{\eps_k}})\chi_{\Omega_2} + \BBB (\id + {\eps_k}^{1/4})\chi_{\Omega_3} \EEE$$for $\alpha \in \R^2$. Then for $b \in \R^2$ possible alternatives are e.g. 
\BBB
\begin{align*}
(1)& \ \ P^1_1 = \Omega_1, \ P^1_2 = \Omega_2, \ P^1_3 = \Omega_3 \ \ \text{with} \ \ T^1_k = y_k \ \text{on} \ \Omega,\\
(2)& \ \ P^2_1 = \Omega_1 \cup \Omega_2, \ P^2_2 = \Omega_3 \ \ \text{with} \ \  T^2_k = \id \chi_{\Omega_1\cup \Omega_2} + (\id + {\eps_k}^{1/4} -b \sqrt{\eps_k})\chi_{\Omega_3}.
\end{align*}
\EEE Letting $u^i_{k} = {\eps_k}^{-\frac{1}{2}}\big(y_{k} - T_k^i\big)$ for $i=1,2$ we obtain in the limit $\eps_k \to 0$  \BBB the unique rigid motion $T=\id$  \EEE and the different configurations \BBB 
\begin{align*}
&(1) \ \  u^1 = 0, \ \ \ \  \ \ \ \ \ \ \ \ \ \ \ \ \ \ \ \ \ \ \ \ \  \ \ \ \ \ \ \ \ \ P^1_1 = \Omega_1, \  P^1_2 = \Omega_2, \ P^1_3 = \Omega_3,\\
 &(2) \ \ u^2 =  0 \cdot \chi_{\Omega_1} + \alpha\chi_{\Omega_2} + b\chi_{\Omega_3}, \ \ \ \ \ \ \ P^2_1 = \Omega_1 \cup \Omega_2, \ P^2_2 = \Omega_3. 
  \end{align*}\EEE}
\end{example}
\vspace{-0.8cm}
We now introduce a special subclass of partitions in which uniqueness will be guaranteed. The above example already shows that different partitions are not equivalent in the sense that they may contain a different `amount of information'. Note that on the various elements of the partition the configuration $u$ is defined separately and the different pieces of the domain are not `aware of each other'. In particular, \BBB the difference of the traces of $u$ on $\partial^* P_i \cap \partial^* P_j$, $i\neq j$, \EEE does not have any physically reasonable interpretation. On the contrary, in example \BBB (2) \EEE where we did not split up \BBB $\Omega_1 \cup \Omega_2$, \EEE we gain the jump height on \BBB $\partial \Omega_1 \cap \partial \Omega_2$ \EEE as an additional information. The observation that coarser partitions provide more information about the behavior at the jump set motivates the definition of the \textit{coarsest partition}. 

\begin{definition}\label{rig-def: ad,coar} 
{\normalfont
Let $(y_k)_k$ be a given (sub-)sequence as in Theorem \ref{rig-th: comp1}.
\begin{itemize}
\item[(i)] We say a partition ${\cal P}$ of $\Omega$ is \textit{admissible} for $(y_k)_k$, and write ${\cal P} \in {\cal Z}_P((y_k)_k)$, if there exist $u,T$ such that $(u, {\cal P}, T) \in {\cal D}_\infty$ and \BBB $y_k \to (u,{\cal P},T)$. \EEE 
\item[(ii)] We say a piecewise rigid motion $T$ is \textit{admissible} for \BBB $(y_k)_k$, \EEE and write  $T \in {\cal Z}_T((y_k)_k)$, if there exist $u, \BBB {\cal P} \EEE $ such that $(u, {\cal P}, T) \in {\cal D}_\infty$ and \BBB $y_k \to (u,{\cal P},T)$\EEE .
\item[(iii)] We say a configuration $u$ is \textit{admissible} for $(y_k)_k$ and ${\cal P}$, and write $u \in {\cal Z}_u((y_k)_k,{\cal P})$, if there exists $T$ such that $(u, {\cal P}, T) \in {\cal D}_\infty$ and \BBB $y_k  \to (u,{\cal P},T)$.\EEE
\item[(iv)] We say a partition  ${\cal P}$ of $\Omega$ is a \textit{coarsest partition} for $(y_k)_k$ if the following holds: The partition is admissible, i.e. ${\cal P} \in {\cal Z}_P((y_k)_k)$. Moreover, for all admissible $u \in {\cal Z}_u((y_k)_k, {\cal P})$ and \BBB all corresponding triples $(u_k,{\cal P}_k,T_k) \in {\cal D}$ satisfying \eqref{rig-eq: comp2}-\eqref{rig-eq: comp3} \EEE the mappings $T_k = \sum_j(R_j^k \cdot + b_j^k)\chi_{P^k_j}$ fulfill
\begin{align}\label{rig-eq: toinfty}
\frac{|R^k_{i} - R^k_{j}| +  |b^k_{i} - b^k_{j}|}{\sqrt{\eps_k}} \to \infty 
\end{align} 
for all $i,j \in \N$,  $i \neq j$ and $k \to \infty$.

\end{itemize}

}
\end{definition}

In Lemma  \ref{rig-lemma: comp2.1} below we find an equivalent characterization of coarsest partitions being the maximal elements of the partial order on the sets of admissible partitions which is induced by subordination. Loosely speaking, the above definition particularly implies that given a coarsest partition  a region of the domain is partitioned into different sets $(P_j)_j$ if and only if the \BBB (scaled) jump height $\eps_k^{-1/2}[y_k]$ \EEE on $(\partial^* P_j)_j$ tends to infinity \BBB (cf. Example \ref{ex}). \EEE

Recall the definition of the  piecewise infinitesimal rigid motions ${\cal A}({\cal P})$ below \eqref{rig-eq: defA}. We now obtain a unique characterization of the limiting configuration up to piecewise infinitesimal rigid motions.

\begin{theorem}\label{rig-th: comp2}
Let $\eps_k \to 0$ be given. \BBB Let $(y_k)_k \subset SBV_M(\Omega)$ be a sequence \EEE for which the assertion of Theorem \ref{rig-th: comp1} holds. Then we have the following: 

\begin{itemize}
\item[(i)] There is a unique $T \in {\cal Z}_T((\BBB y_{k})_k \EEE )$.
\item[(ii)] There is a unique coarsest partition $\bar{\cal P}$ of $\Omega$. 
\item[(iii)] Given some $u \in {\cal Z}_u((y_{k})_k, \bar{\cal P})$ all admissible limiting configurations are of the form $u + \nabla T{\cal A}(\bar{\cal P})$, i.e. the limiting configuration is determined uniquely up to piecewise infinitesimal rigid motions. 
\end{itemize}
\end{theorem}

\BBB Going back to Example \ref{ex}, we observe that $T = \id$ is uniquely given and that the partition in (2) is the coarsest partition. The non-uniqueness in (iii) is a consequence of the fact that the nonlinear energy is invariant under rigid motions (see also (2) in Example \ref{ex} for $b \in \R^2$). \EEE

%----------------------------------------------------------------------------------------------------------------------
\subsection{The limiting linearized model and $\Gamma$-convergence}\label{rig-sec: sub, main gamma}
%-----------------------------------------------------------------------------------------------------------------------

We now introduce the limiting linearized model, discuss  its properties and show that it can be identified as the $\Gamma$-limit of the nonlinear energies $E_\eps$. Let $Q = D^2W(\Id)$ be the Hessian of the stored energy density $W$ at the identity. Define $E : {\cal D}_\infty \to [0,\infty)$ by
\begin{align}\label{rig-eq: Griffith en-lim}
E(u,{\cal P},T) =  \int_\Omega \frac{1}{2} Q(e(\nabla T^T \nabla u))  + {\cal H}^1\Big(J_u \setminus \bigcup\nolimits_j\partial^* P_j\Big) + {\cal H}^1\Big(\bigcup\nolimits_j \partial^* P_j \cap \Omega\Big),
\end{align} 
where as before ${\cal P} = (P_j)_j$. Recall that a triple of the limiting model consists of a partition, a corresponding piecewise rigid motion and a displacement field. \BBB We emphasize that in contrast to the nonlinear model (see \eqref{rig-eq: SBVfirstdef} and Remark \ref{rem:M}) there are no restrictive bounds on the functions $u$ and their derivatives\EEE .  

The surface energy of $E$ has two parts. Similarly as discussed in Section \ref{rig-sec: sub, seg}, on the right we have the \textit{segmentation energy} which is necessary to disconnect the components of the body. Moreover, on the left we have the \textit{inner crack energy} associated to the discontinuity set of the displacement field in each part of the material \BBB (see also Remark \ref{rem:new}(ii) below). \EEE Whereas the first two terms of the functional typically appear in the study of linearized Griffith energies, the segmentation energy is a characteristic feature of our general model where the analysis is not restricted to a linearization around a fixed rigid motion. 

We now present our main $\Gamma$-convergence result.  Recall \BBB Definition \ref{def:conv}. \EEE

\begin{theorem}\label{rig-th: gammaconv}
Let $\Omega \subset \R^2$ open, bounded with Lipschitz boundary. Let $M> 0$ and $\eps_k \to 0$. Then \BBB $E_{\eps_k}$ converges to $E$ in the sense of $\Gamma$-convergence, \EEE i.e.
\begin{itemize}
\item[(i)] $\Gamma-\liminf$ inequality: For all $(u,{\cal P},T) \in {\cal D}_\infty$ and for all sequences $(y_k)_k \subset SBV_M(\Omega)$ \BBB with $y_k \to (u,{\cal P},T)$ \EEE we have 
$$\liminf_{k \to \infty} E_{\eps_k}(y_k) \ge E(u,{\cal P},T).$$
\item[(ii)] Existence of recovery sequences: For every $(u,{\cal P},T) \in {\cal D}_\infty$  with $u \in L^2(\Omega)$ we find a sequence $(y_k)_k \subset SBV_M(\Omega)$ \BBB such that $y_k \to (u,{\cal P},T)$ and \EEE
$$\lim_{k \to \infty} E_{\eps_k}(y_k) = E(u,{\cal P},T).$$
\end{itemize}
\end{theorem} 

\vspace{-0.6cm}
\begin{rem}\label{rem:new}
{\normalfont
(i) \BBB The limiting model could equivalently be formulated with $v = \nabla T^T u$ in place of the displacement field $u$. (Accordingly, replace $u_k$ by $v_k = \nabla T_k^T u_k$ in Theorem \ref{rig-th: comp1}). This alternative notation simplifies the description of the elastic energy in \eqref{rig-eq: Griffith en-lim}, but does not account for the fact that the linearization was possibly performed around different rigid motions. 

(ii) Using  the local structure of Caccioppoli partitions (see Theorem \ref{th: local structure} and recall \eqref{eq: essential boundary}) the limiting energy can equivalently be written as
$$\sum\nolimits_j \Big( \int_{P_j} \frac{1}{2} Q(e(R_j^T \nabla u))  + \mathcal{H}^1(J_u \cap (P_j)^1) + \frac{1}{2}\mathcal{H}^1(\partial^* P_j \cap \Omega) \Big).$$ \EEE

(iii) For configurations $(u,\bar{\cal P},T)$ defined in terms of the coarsest partition $\bar{\cal P}$ there is an additional interpretation for the  crack opening of the sequence of deformations $y_\eps$: (1) The jumps on $\bigcup_j\partial^* P_j$ are associated to jump heights $\gg \sqrt{\eps}$ and (2) the jump heights corresponding to the inner crack energy are of the order $\sqrt{\eps}$, \BBB which illustrates the multiscale nature of the model. \EEE In fact, (1) follows from \eqref{rig-eq: toinfty} and (2) is a consequence of \eqref{rig-eq: comp2}(i).

(iv) On a component $P_j$ of $\bar{\cal P}$ the body may still be disconnected by the jump set \BBB $(P_j)^1 \cap J_u$ \EEE forming a finer partition of the specimen. However, in contrast to the boundary of  $\bar{\cal P}$ the jump heights  have a meaningful physical interpretation. 

\BBB (v)  In general, the partition induced by the \emph{macroscopic jumps} (represented by $J_T$) is  coarser than $\bar{\cal P}$, i.e. $\mathcal{H}^1(\bigcup_j\partial^* P_j \setminus (\partial \Omega \cup J_T))>0$, cf. Example \ref{ex}. \EEE

}
\end{rem}

As a direct consequence of Theorem \ref{rig-th: gammaconv} we get that the $\Gamma$-limit \BBB of the same functionals $E_{\eps_k}$ with respect to the much weaker notion of $L^1$-convergence of the unrescaled deformations $y_k$ is given by the segmentation energy. \EEE

\begin{corollary}\label{rig-cor: gamma}
Let $\Omega \subset \R^2$ open, bounded with Lipschitz boundary. Let $M> 0$ and $\eps_k \to 0$. Then  $E_{\eps_k}$   $\Gamma$-converge to $E_{\rm seg}$   with respect to the  $L^1(\Omega)$-convergence, where
$$E_{\rm seg}(y) = \begin{cases} {\cal H}^1\big( \BBB J_T \EEE \big)& y = T \in {\cal R}({\cal P}) \ \text{ for a Caccioppoli partition } {\cal P}, \\ + \infty & \text{else.}\end{cases} $$
\end{corollary}

\BBB Note that the segmentation energy in Corollary \ref{rig-cor: gamma} differs from the one in \eqref{rig-eq: Griffith en-lim}, see Remark \ref{rem:new}(v). 
 
 \BBB
 
 %----------------------------------------------------------------------------------------------------------------------
\subsection{Application: External loads}\label{rig-sec: sub, cleavage2}
%-----------------------------------------------------------------------------------------------------------------------

For the investigation of minimization problems associated to $E_\eps$ it is interesting to take external loads into account. In the context of brittle materials, however, the incorporation of body forces is a delicate problem. Indeed, assumptions on the class of admissible loads have to ensure that no part of the body is broken apart and sent to infinity, which clearly excludes the case of a constant body force (see \cite[Remark 3.1]{DalMaso-Francfort-Toader:2005}).  To avoid the occurrence of such phenomena, it is natural to assume a uniform  $L^\infty$ bound on the admissible functions, see e.g. \cite{DFT2, DalMaso-Lazzaroni:2010}. Unfortunately, this is not expedient  in our setting since the bound $\Vert y \Vert_\infty \le M$ is futile after passage to rescaled configurations   in Theorem \ref{rig-th: comp1}.

Let us mention that in \cite{DalMaso-Francfort-Toader:2005}  body forces were indispensable to ensure reasonable compactness properties for sequences and to guarantee existence of minimizers. In our setting, however, due to the subtraction of suitable rigid motions on a partition of the domain (cf. \eqref{rig-eq: comp2}(i)) we obtain a compactness result without the necessity of  additional loading terms.

We fix a sequence $\eps_k \to 0$ and consider the following prototype problem $F_{\eps_k} : SBV_M(\Omega) \to [0,\infty)$ with
\begin{align}\label{ennew}
F_{\eps_k}(y) = E_{\eps_k}(y) + \frac{\lambda}{\eps_k} \Vert y - f_k \Vert^2_{L^2(\Omega)},
\end{align}
where  $\lambda >0$ and $(f_k)_k \subset SBV_M(\Omega)$ a sequence  with $\sup_k E_{\eps_k}(f_k) < \infty$. An expansion yields the constant $ \lambda\eps_k^{-1} \int_\Omega|f_k|^2$, the \emph{external load} $-2\lambda\eps_k^{-1} \int_\Omega f_k \cdot y$ and the term  $\lambda \eps_k^{-1} \int_\Omega |y|^2$. The latter can be interpreted as an \emph{ artificial confining potential}, which prevents parts of the body from being sent to infinity.   

We assume that there is  a triple $(g,{\cal P}_g,T_g) \in \mathcal{D}_\infty$ such that $f_k \to (g,{\cal P}_g,T_g)$ in the sense of Definition \ref{def:conv} and that for the associated triples $(g_k,{\cal P}_k^g,T_k^g) \in {\cal D}$   satisfying  \eqref{rig-eq: comp2}-\eqref{rig-eq: comp3} we have $\eps_k^{-1/2}(f_k - T_k^g) \to g$ in $L^2(\Omega)$ and ${\cal P}_k^g = {\cal P}_g$ for all $k \in \N$. (Note that up to a subsequence the convergence in the sense of Definition \ref{def:conv} is already guaranteed by Theorem \ref{rig-th: comp1}.) Moreover, we suppose that ${\cal P}_g$ is the coarsest partition given by Theorem \ref{rig-th: comp2}(ii) and write ${\cal P}_g = (P^g_j)_j$. 
By ${\cal C}_g \subset {\cal D}_\infty$ we denote the set of triples  $(u,{\cal P},T) \in {\cal D}_\infty$ with $T=T_g$  and the property that ${\cal P}_g$ is coarser than ${\cal P}$, i.e. for each $P_j$ there exists $P^g_i$ with $|P_j \setminus P_i^g|= 0$.

\begin{lemma}\label{lemmanew}
Let $(y_k)_k \subset SBV_M(\Omega)$ be a sequence with $F_{\eps_k}(y_{k}) \le C$  and $y_k \to (u,{\cal P}, T) \in {\cal D}_\infty$ in the sense of Definition \ref{def:conv}. Then $(u,{\cal P},T) \in {\cal C}_g$.
\end{lemma}

Recalling \eqref{rig-eq: Griffith en-lim} we introduce the limiting energy    $F_g : {\cal D}_\infty \to [0,\infty]$ by
\begin{align*}
F_g(u,{\cal P},T) = \begin{cases} E(u,{\cal P},T) +  \min_{v \in u + \nabla T \mathcal{A}(\mathcal{P})} \lambda\Vert v- g \Vert^2_{L^2(\Omega)} & \text{if }  (u,{\cal P},T) \in {\cal C}_g,\\ + \infty &\text{else}, \end{cases}
\end{align*}
where $\mathcal{A}(\mathcal{P})$ as defined below \eqref{rig-eq: defA}. Similarly as the functional in \eqref{rig-eq: Griffith en-lim}, $F_g$ is invariant under infinitesimal rigid motions on the components of the partition ${\cal P}$. However, the additional term on the right induces a symmetry breaking and there is exactly one distinguished configuration $u^*$ in the class $u + \nabla T \mathcal{A}(\mathcal{P})$ (cf. (iii) in Theorem \ref{rig-th: comp2}) which satisfies $\min_{v \in u + \nabla T \mathcal{A}(\mathcal{P})}\Vert v- g \Vert^2_{L^2(\Omega)} = $ $  \Vert u^*- g \Vert^2_{L^2(\Omega)}$. We close this section with a corresponding $\Gamma$-convergence result.

\begin{theorem}\label{rig-th: gammaconv2}
Let $\Omega \subset \R^2$ open, bounded  with Lipschitz boundary. Let $M> 0$, $\eps_k \to 0$  and $(f_k)_k, g$ as above. Then  $F_{\eps_k}$ converges to $F_g$ in the sense of $\Gamma$-convergence. (Replace $E_{\eps_k}$ by $F_{\eps_k}$ and $E$ by  $F_g$ in (i),(ii) of Theorem \ref{rig-th: gammaconv}.) Moreover,  we  have 
$$\lim_{k \to \infty} \ \  \inf_{y \in SBV_M(\Omega)} F_{\eps_k}(y) = \min_{(u,{\cal P},T) \in {\cal D}_\infty} F_g(u,{\cal P},T)$$ and (almost) minimzers of $F_{\eps_k}$ converge (up to subsequences) to minimizers of $F_g$ in the sense of Definition \ref{def:conv}.
\end{theorem}

\EEE

%----------------------------------------------------------------------------------------------------------------------
\subsection{\BBB Application: Cleavage laws \EEE }\label{rig-sec: sub, cleavage}
%-----------------------------------------------------------------------------------------------------------------------

In fracture mechanics it is a major challenge to identify critical loads at which
a body fails and to determine the geometry of crack paths that occur in the
fractured regime. As \BBB another \EEE application of the above results we now finally derive such a cleavage law. We consider a special boundary value problem of uniaxial compression/extension. Let $\Omega = (0,l) \times (0,1)$, $\Omega' = (-\eta, l + \eta) \times (0,1)$ for $l>0$, $\eta >0$ and for $a_\eps \in \R$ define
$${\cal A}(a_\eps) := \lbrace y \in SBV_M(\Omega'): y_1 (x) = (1+ a_\eps)x_1 \text{ for } x_1 \le 0 \ \text{or} \ x_1 \ge l \rbrace.$$  
As usual in the theory of $SBV$ functions the boundary values have to be imposed in small neighborhoods of the boundary. In what follows,  the elastic part of the energy \eqref{rig-eq: Griffith en} still only depends on $y|_\Omega$,  whereas the surface energy is given by ${\cal H}^1(J_y)$ with $J_y \subset \Omega'$. In particular, jumps on $\lbrace0,l\rbrace \times (0,1)$ contribute to the energy $E_\eps(y)$. (Also compare a similar discussion before \cite[Theorem 2.2]{FriedrichSchmidt:2014.2}.) The present problem in the framework of continuum fracture mechanics with isotropic surface energies is a slightly simplified model of the problem considered in \cite{FriedrichSchmidt:2011, FriedrichSchmidt:2014.2}.  

As a preparation, define $\alpha$ such that $\inf \lbrace Q(F): \e_1^T F \e_1 = 1\rbrace = \alpha$ and observe $\inf \lbrace Q(F): \e_1^T F \e_1 = a\rbrace = \alpha a^2$ for all $a \in \R$. Moreover, let $F^a \in \R^{2 \times 2}_{\rm sym}$ be the unique matrix such that $\e_1^T F^a \e_1 = a$ and $Q(F^a) = \inf \lbrace Q(F): \e_1^T F \e_1 = a\rbrace = \alpha a^2$. 

We recall that the proof of the cleavage laws in \cite{FriedrichSchmidt:2011, FriedrichSchmidt:2014.1, Mora:2010}  fundamentally relied on the application of certain slicing techniques which were not suitable to treat the case of compression. Having general compactness and $\Gamma$-convergence results we can now  complete the picture about cleavage laws by extending the results to the case of uniaxial compression.

\begin{theorem}\label{rig-th: cleavage-cont}
Suppose $a_\eps/\sqrt{\eps} \to a \in [-\infty,\infty]$. The limiting minimal energy is given by
\begin{align}\label{rig-eq: cleavage en}
\lim_{\eps \to 0} \inf \lbrace E_\eps(y): y \in {\cal A}(a_\eps)\rbrace = \min\Big\{  \frac{1}{2} \alpha l a^2,1\Big\}.
\end{align}
Let $a_{\rm crit}:= \sqrt{\frac{2\alpha}{l}}$.  For every sequence $(y_\eps)_\eps$ of almost minimizers, up to passing to subsequences, we get $\eps^{-1/2}(y_\eps(x) - x) \to u(x)$ for a.e. $x \in \Omega$, where
\begin{itemize}
\item[(i)] if $|a| < a_{\rm crit}$, $u(x) = (0,s) + F^a x$ for $s\in \R$,
\item[(ii)] if $|a| > a_{\rm crit}$,  $u(x) = \begin{cases} (0,s) & x_1 < p, \\ (l a,t) & x_1 > p, \end{cases}$ for $s,t \in \R$, $p\in (0,l)$.
\end{itemize}
\end{theorem} 

\BBB Let us emphasize that the cleavage law is derived   for a special geometry of $\Omega$ by solving a static, global minimization problem similarly as in  \cite{Braides-Lew-Ortiz:06, FriedrichSchmidt:2014.1, Mora:2010}. An accurate prediction of crack propagation under tensile loading is beyond the scope of the present contribution. \EEE

%----------------------------------------------------------------------------------------------------------------------
\section{Preliminaries}\label{rig-sec: pre}
%-----------------------------------------------------------------------------------------------------------------------

In this section  we collect the definitions as well as basic properties of $SBV$ and $SBD$ functions and state the rigidity estimates which are necessary for the derivation of our main compactness result. 

%----------------------------------------------------------------------------------------------------------------------
\subsection{(G)SBV and (G)SBD functions}\label{rig-sec: sub, sbv}
%-----------------------------------------------------------------------------------------------------------------------

Let $\Omega \subset \R^d$ open, bounded with Lipschitz boundary. Recall that the space $SBV(\Omega, \R^d)$, abbreviated as $SBV(\Omega)$ hereafter,  of \emph{special functions of bounded variation} consists of functions $y \in L^1(\Omega, \R^d)$ whose distributional derivative $Dy$ is a finite Radon measure, which splits into an absolutely continuous part with density $\nabla y$ with respect to Lebesgue measure and a singular part $D^s y$ whose Cantor part vanishes and thus is of the form 
$$ D^s y = [y] \otimes \nu_y {\cal H}^{d-1} \lfloor J_y. $$
Here ${\cal H}^{d-1}$ denotes the $(d-1)$-dimensional Hausdorff measure, $J_y$ (the `crack path') is an ${\cal H}^{d-1}$-rectifiable set in $\Omega$, $\nu_y$ is a normal of $J_y$ and $[y] = y^+ - y^-$ (the `crack opening') with $y^{\pm}$ being the one-sided limits of $y$ at $J_y$. If in addition $\nabla y \in L^2(\Omega)$ and ${\cal H}^{d-1}(J_y) < \infty$, we write $y \in SBV^2(\Omega)$. See \cite{Ambrosio-Fusco-Pallara:2000} for the basic properties of this function space. 

Likewise, we say that a function $y \in L^1(\Omega, \R^d)$ is a \emph{special  function of bounded deformation} if the symmetrized distributional derivative $Ey := \frac{(Dy)^T + Dy}{2}$ is a finite \BBB $\R^{d \times d}_{\rm sym}$-valued \EEE Radon measure with vanishing Cantor part. It can be decomposed as 
\begin{align}\label{rig-eq: symmeas}
 Ey = e(\nabla y) {\cal L}^d  + E^s y = e(\nabla y) {\cal L}^d + [y] \odot \nu_y {\cal H}^{d-1}|_{J_y},
 \end{align}
where $e(\nabla y)$ is the absolutely continuous part of $Ey$ with respect to the Lebesgue measure ${\cal L}^d$, $[y]$, $\nu_y$, $J_y$ as before and $a \odot b = \frac{1}{2}(a \otimes b + b \otimes a)$. For basic properties of this function space we refer to \cite
{Ambrosio-Coscia-Dal Maso:1997,  Bellettini-Coscia-DalMaso:98}.

To treat variational problems as considered in Section \ref{rig-sec: main} (see in particular \eqref{rig-eq: Griffith en}) the spaces $SBV(\Omega)$ and $SBD(\Omega)$ are not adequate due to the lacking $L^\infty$-bound being essential in the compactness theorems. To overcome this difficulty the space of $GSBV(\Omega)$ was introduced consisting of all ${\cal L}^d$-measurable functions $y: \Omega \to \R^d$ such that for every $\phi \in C^1(\R^d)$ with the support of $\nabla \phi$ compact, the composition $\phi \circ y $ belongs to $SBV_{\rm loc}(\Omega)$ (see \cite{DeGiorgi-Ambrosio:1988}). In this new setting one may obtain  a more general compactness result (see \cite[Theorem 4.36]{Ambrosio-Fusco-Pallara:2000}). Unfortunately, this approach cannot be pursued in the framework of $SBD$ functions as for a function $y \in SBD(\Omega)$ \BBB the composition \EEE $\phi \circ y$ typically does not lie in $SBD(\Omega)$. In \cite{DalMaso:13}, Dal Maso suggested another approach which is based on certain properties of one-dimensional slices. 

First we have to introduce some notation. For every $\xi \in \R^d \setminus \lbrace 0 \rbrace$, for every $s \in \R^d$ and for every $B \subset \Omega$ we let 
\begin{align}\label{eq: slicing 1}
B^{\xi,s} = \lbrace t \in \R: s + t\xi \in B\rbrace.
\end{align}
Furthermore, define the hyperplane $\Pi^\xi = \lbrace x \in \R^d: x \cdot \xi = 0\rbrace$. Moreover, for every function $y: B \to \R^d$ we introduce the function $y^{\xi,s} : B^{\xi,s} \to \R^d$ by 
\begin{align}\label{eq: slicing 2}
y^{\xi,s}(t) = y(s + t\xi)
\end{align}
 and  $\hat{y}^{\xi,s} : B^{\xi,s} \to \R$ by $\hat{y}^{\xi,s}(t) = y(s + t\xi) \cdot \xi$. If $\hat{y}^{\xi,s} \in SBV(B^{\xi,s},\R)$ and $J_{\hat{y}^{\xi,s}}$ denotes \BBB the \EEE \textit{approximate jump set}, we define 
$$J^1_{\hat{y}^{\xi,s}} := \lbrace t \in J_{\hat{y}^{\xi,s}}: |[\hat{y}^{\xi,s}]( t)| \ge 1 \rbrace.$$
The space $GSBD(\Omega,\R^d)$ of \emph{generalized functions of bounded deformation} is the space of all ${\cal L}^d$-measurable functions $y: \Omega \to \R^d$ with the following property: There exists a nonnegative bounded Radon measure $\lambda$ on $\Omega$ such that for all $\xi \in S^{d-1}:=\lbrace x \in \R^d: |x|=1 \rbrace$ we have that for ${\cal H}^{d-1}$-a.e. $s \in \Pi^\xi$ the function $\hat{y}^{\xi,s} = y^{\xi,s} \cdot \xi$ belongs to $SBV_{\rm loc}(\Omega^{\xi,s})$ and
$$\int_{\Pi^\xi} \Big( |D\hat{y}^{\xi,s}|(B^{\xi,s} \setminus J^1_{\hat{y}^{\xi,s}}) + {\cal H}^0(B^{\xi,s} \cap J^1_{\hat{y}^{\xi,s}})\Big)\, d{\cal H}^{d-1}(s) \le \lambda(B)$$
for all Borel sets $B \subset \Omega$.

We refer to \cite{DalMaso:13} for basic properties of this space. In particular, for later reference we now recall fundamental slicing, compactness and approximation results. We first briefly state the main slicing properties of $GSBD$ functions (see \cite[Section 8,9]{DalMaso:13}.) Recall definitions \eqref{eq: slicing 1} and \eqref{eq: slicing 2} and let $J_y^\xi = \lbrace x \in J_y: [y](x) \cdot \xi \neq 0 \rbrace$. 

\begin{theorem}\label{clea-th: slic}
Let $y \in GSBD(\Omega)$. For all $\xi \in S^{d-1}$ and ${\cal H}^{d-1}$-a.e.\ $s$ in $\Pi^\xi = \lbrace x: x\cdot \xi = 0\rbrace$ we have $J_{\hat{y}^{\xi,s}} = (J^\xi_y)^{\xi,s}$ and
 \begin{align*}
 \int_{\Pi^\xi} \# J_{\hat{y}^{\xi,s}} \, d{\cal H}^{d-1}(s) = \int_{J^\xi_y} |\nu_y \cdot \xi| \, d{\cal H}^{d-1}.
 \end{align*}
Moreover, the approximate symmetrized gradient $e(\nabla y)$ exists in the sense of \cite[(9.1)]{DalMaso:13}, satisfies $e(\nabla y) \in L^1(\Omega, \R_{\rm sym}^{d \times d})$ and for all $\xi \in S^{d-1}$ and ${\cal H}^{d-1}$-a.e.\ $s$ in $\Pi^\xi$ we have
\begin{align*}
\xi^T e(\nabla y(s+ t\xi))\xi =  (\hat{y}^{\xi,s})'(t) \ \text{ for a.e. } t \in {\Omega}^{\xi,s}.
 \end{align*}
\end{theorem}
 
If in addition $e(\nabla y) \in L^2(\Omega)$ and ${\cal H}^{d-1}(J_y) < \infty$, we write $y \in GSBD^2(\Omega)$.  Similar properties for $SBV$ functions may be found in \cite[Section 3.11]{Ambrosio-Fusco-Pallara:2000}. We now state a general compactness result in $GSBD$ proved in \cite[Theorem 11.3]{DalMaso:13} which we slightly adapt for our purposes.

\begin{theorem}\label{rig-th: GSBD comp}
Let $(y_k)_k$ be a sequence in \BBB $GSBD^2(\Omega)$. \EEE Suppose that there exist a constant $M>0$ and an increasing continuous functions $\psi:[0,\infty) \to [0,\infty)$ with $\lim_{ \BBB t \EEE \to \infty} \psi( \BBB t \EEE ) = + \infty$ such that 
$$\int_{\Omega} \psi(|y_k|) + \int_{\Omega} |e(\nabla y_k)|^2 + {\cal H}^{d-1}(J_{y_k}) \le M $$
for every $k \in \N$. Then there exist a subsequence, still denoted by $(y_k)_k$, and a function $y \in GSBD^2(\Omega
)$ such that
\begin{align}\label{rig-eq: convergence sense}
\begin{split}
& y_k \to y \ \ \ \text{pointwise a.e. in} \ \ \ \Omega,\\ 
&e (\nabla y_k) \rightharpoonup  e (\nabla y)  \ \ \text{ weakly in} \ L^2(\Omega,\R^{d\times d}_{\rm sym}),\\
& \liminf_{k \to \infty} {\cal H}^{d-1}(J_{y_k}) \ge {\cal H}^{d-1}(J_y).
\end{split}
\end{align}
\end{theorem}

The lower semicontinuity result for the jump set can be generalized considering one-dimensional slices. Define  $\theta_\sigma: [0,\infty) \to [0,1]$ by $\theta_\sigma(t) = \min \lbrace \frac{t}{\sigma},1 \rbrace$ for  $\sigma > 0$ and additionally $\theta_0 \equiv 1$. Let
\begin{align}\label{rig-eq: lemma2**}
\hat{\mu}^{\sigma,\xi}_y(B) := \int_{\Pi^\xi} \int_{B^{\xi,s} \cap J_{\hat{y}^{\xi,s}}} \theta_\sigma(|[\hat{y}^{\xi,s}](t) | )\, d{\cal H}^0(t) \, d{\cal H}^{d-1}(s)
\end{align}
for all Borel sets $B \subset \Omega$.

\begin{lemma}\label{rig-lemma: lemma}
Let $(y_k)_k$ be a sequence in \BBB $GSBD^2(\Omega)$ \EEE converging to a function $y \in \BBB GSBD^2(\Omega) \EEE $ in the sense of \eqref{rig-eq: convergence sense}. Then
\begin{align*}
\hat{\mu}^{\sigma,\xi}_y(U) \le \liminf_{k \to \infty} \hat{\mu}^{\sigma,\xi}_{y_k}(U) 
\end{align*}
for all $\sigma \ge 0$, every $\xi \in S^{d-1}$ and for all open sets $U \subset \Omega$.  
\end{lemma} 

\Proof As $y_k \to y$ in the sense of \eqref{rig-eq: convergence sense}, we may assume that  $ (\hat{y}_k)^{\xi,s}  \to \hat{y}^{\xi,s}$ in $GSBV(U^{\xi,s})$ for ${\cal H}^{d-1}$-a.e. $s \in U^\xi := \lbrace s \in \Pi^\xi: U^{\xi,s} \neq \emptyset\rbrace$. This is one of the essential steps in the proof of Theorem \ref{rig-th: GSBD comp}  (cf. \cite[Theorem 11.3]{DalMaso:13} or \cite[Theorem 1.1]{Bellettini-Coscia-DalMaso:98} for an elaborated proof in the $SBD$-setting). The desired claim now follows from the corresponding lower semicontinuity result for $GSBV$ functions (see e.g. \cite[Theorem 4.36]{Ambrosio-Fusco-Pallara:2000}) and Fatou's lemma. \eop

We briefly note that using the area formula (see e.g. \cite[Theorem 2.71]{Ambrosio-Fusco-Pallara:2000})) and fine properties of $GSBD$ functions (see \cite{DalMaso:13}), $\hat{\mu}^{\sigma,\xi}_y(B)$ can be written equivalently as
\begin{align}\label{rig-eq: lemma2}
\hat{\mu}^{\sigma,\xi}_y(B) = \int_{\BBB J^\xi_y \EEE \cap B} \theta_\sigma(|[y] \cdot \xi|) |\nu_y \cdot \xi| \, d{\cal H}^{d-1}
\end{align}
for all \BBB $\sigma \ge 0$, \EEE for all $\xi \in S^{d-1}$ and all Borel sets $B \subset \Omega$ (see also \cite[Remark 9.3]{DalMaso:13}). Finally,  we recall a density result in  $GSBD$ (see \cite{Iurlano:13}).

\begin{theorem}\label{rig-th: cortesani2}
Let $y \in GSBD^2(\Omega) \cap L^2(\Omega)$. Then there exists a sequence $y_k \in SBV^2(\Omega)$ such that each $J_{y_k}$ is contained in the union of  a finite number of closed connected pieces of $C^{1}$-hypersurfaces, each $y_k$ belongs to
$W^{1,\infty}(\Omega \setminus \overline{J_{y_k}},\R^d)$ and the following properties hold:
\begin{align*}
(i) & \ \ \Vert y_k - y \Vert_{L^2(\Omega)} \to 0,\\
(ii) & \ \ \Vert e(\nabla y_k) - e(\nabla y) \Vert_{L^2(\Omega)} \to 0,\\
(iii) &  \ \ {\cal H}^{d-1}(J_{y_k} \BBB \triangle \EEE J_y) \to  0.
\end{align*}
\end{theorem}

%----------------------------------------------------------------------------------------------------------------------
\subsection{Caccioppoli partitions}\label{rig-sec: sub, cacciop}
%-----------------------------------------------------------------------------------------------------------------------

 Let $\Omega \subset \R^d$ open and $E \subset \Omega$ measurable. For \BBB $t \in [0,1]$ we define the points of density $t$ by 
\begin{align}\label{eq: essential boundary}
E^t =   \left\{ x \in \R^d: \lim\nolimits_{\varrho \downarrow 0} \frac{|E \cap B_\varrho(x)|}{|B_\varrho(x)|}  = t\right\}
\end{align}
(see \cite[Definition 3.60]{Ambrosio-Fusco-Pallara:2000}). By $\partial^* E = \R^d \setminus (E^0 \cup E^1)$ we denote the \emph{essential boundary} of $E$ and $\mathcal{H}^1(\partial^*E \cap \Omega)$ denotes the \emph{perimeter} of $E$ in $\Omega$ (cf. \cite[(3.62)]{Ambrosio-Fusco-Pallara:2000})\EEE .

We say a partition ${\cal P} = (P_j)_{j\in\N}$ of $\Omega$ is a \textit{Caccioppoli partition} of $\Omega$  if $\sum_j \mathcal{H}^1(\partial^*P_j) < + \infty$. We say a partition is \textit{ordered} if $|P_i| \ge |P_j|$ for $i \le j$. In the whole paper we will always tacitly assume that partitions are ordered. Given a rectifiable  set $S$ we say that a Caccioppoli partition is \textit{subordinated} to $S$ if (up to an ${\cal H}^{d-1}$-negligible set) the essential boundary $\partial^* P_j$ of $P_j$ is contained in $S$ for every $j \in \N$.  

The  local structure of Caccioppoli partitions can be characterized as follows (see \cite[Theorem 4.17]{Ambrosio-Fusco-Pallara:2000}).
\begin{theorem}\label{th: local structure}
Let $(P_j)_j$ be a Caccioppoli partition of $\Omega$. Then 
$$\bigcup\nolimits_j (P_j)^1 \cup \bigcup\nolimits_{i \neq j} \partial^* P_i \cap \partial^* P_j$$
contains ${\cal H}^{d-1}$-almost all of $\Omega$, \BBB where $(P_j)^1$ as defined in  \eqref{eq: essential boundary}. \EEE
\end{theorem}
 Essentially, the  theorem states that ${\cal H}^{d-1}$-a.e. point of $\Omega$ either belongs to exactly one element of the partition or to the intersection of exactly two sets $\partial^* P_i$, $\partial^* P_j$. We now state a compactness result for ordered Caccioppoli partitions (see \cite[Theorem 4.19, Remark 4.20]{Ambrosio-Fusco-Pallara:2000}).

\begin{theorem}\label{th: comp cacciop}
Let $\Omega \subset \R^d$ open, bounded with Lipschitz boundary.  Let ${\cal P}_i = (P_{j,i})_j$, $i \in \N$, be a sequence of ordered Caccioppoli partitions of $\Omega$ fulfilling $\sup_i \sum_j \mathcal{H}^{d-1}(\partial^* P_{j,i}) < \infty$. Then there exists a Caccioppoli partition ${\cal P} = (P_j)_j$ and a not relabeled subsequence such that \BBB $|P_{j,i} \triangle  P_j| \to 0$ \EEE for all $j \in \N$ as $i \to \infty$.
\end{theorem}
We will also use the fact that $|P_{j,i} \triangle  P_j| \to 0$  for all $j \in \N$ is equivalent to $\sum_j|P_{j,i} \triangle  P_j| \to 0$.  Caccioppoli partitions are naturally associated to piecewise constant functions. We say $y: \Omega \to \BBB \R^m \EEE $ is \emph{piecewiese constant in  $\Omega$} if there exists a Caccioppoli partition $(P_j)_j$ of $\Omega$ and a sequence $(t_j)_j \subset \BBB \R^m \EEE $ such that $y = \sum_j t_j \chi_{P_j}$. We close this section with a compactness result for piecewise constant functions (see \cite[Theorem 4.25]{Ambrosio-Fusco-Pallara:2000}).  

\begin{theorem}\label{th: piecewise const}
Let $\Omega \subset \R^d$ open, bounded with Lipschitz boundary. Let $(y_i)_i \subset SBV(\Omega, \BBB \R^m \EEE )$ be a sequence of piecewise constant functions such that $\sup_i (\Vert y_i \Vert_\infty + {\cal H}^{d-1}(J_{y_i})) < \infty$. Then there exists a  not relabeled subsequence converging in measure to a piecewise constant function $y$.
\end{theorem}

\subsection{Rigidity estimates}\label{sec: rig}

In this section we first recall a geometric rigidity result obtained in the framework of nonlinear elasticity  and a piecewise rigidity estimate for brittle materials for the sake of completeness. Afterwards we introduce a quantitative result in $SBD$ adapted for Griffith energies of the form \eqref{rig-eq: Griffith en} which will be the starting point for our analysis.

We begin with the quantitative geometric rigidity result by Friesecke, James, M\"uller \cite{FrieseckeJamesMueller:02} generalizing the classical Liouville theorem.

\begin{theorem}\label{rig-th: geo rig}
Let $\Omega \subset \R^d$ a (connected) Lipschitz domain and $1 < p < \infty$. Then there exists a constant $C = C(\Omega,p)$ such that for any $y \in W^{1,p}(\Omega,\R^d)$ there is a rotation $R \in SO(d)$ such that
\begin{align*}
\left\|\nabla y - R\right\|_{L^p(\Omega)} \leq C \left\|\dist(\nabla y, SO(d))\right\|_{L^p(\Omega)}. 
\end{align*}
\end{theorem}

In the theory of fracture mechanics the problem is more involved as global rigidity can fail if the crack disconnects the body. Chambolle, Giacomini and Ponsiglione \cite{Chambolle-Giacomini-Ponsiglione:2007} have proved the following qualitative result for brittle materials which do not store elastic energy (i.e. $\nabla y \in SO(d)$  a.e. in $\Omega$). 

\begin{theorem}\label{rig-cor: cgp}
Let $y \in SBV(\Omega)$ such that ${\cal H}^{d-1}(J_y) < +\infty$ and $\nabla y \in SO(d)$  a.e. in $\Omega$. Then $y$ is a collection of an at most countable family of rigid deformations, i.e., there exists a Caccioppoli partition ${\cal P} = (P_j)_j$ subordinated to $J_y$ such that
$$y(x) = \sum\nolimits_j (R_j \, x + b_j) \chi_{P_j}(x),$$
where $R_j \in SO(d)$ and $b_j \in \R^d$. 
\end{theorem}

Loosely speaking, the result states that the only way that rigidity may fail is that the body is divided into at most countably many parts each of which subject to a different rigid motion. We briefly note that there is an analogous result in the geometrically linear setting (see \cite[Theorem A.1]{Chambolle-Giacomini-Ponsiglione:2007}): A function $u \in SBD(\Omega)$ with ${\cal H}^{d-1}(J_u) < +\infty$ and $e(\nabla u) = 0$  a.e. in $\Omega$ has the form $u(x) = \sum\nolimits_j (A_j \, x + b_j) \chi_{P_j}(x)$  for $A_j \in \R^{d \times d}_{\rm skew}$ and $b_j \in \R^d$.

We now introduce a quantitative $SBD$-rigidity result which may be seen as a suitable combination of the above estimates and is tailor-made for general Griffith functionals of the form \eqref{rig-eq: Griffith en} where both energy forms are coexistent (see \cite[Theorem 2.1, Remark 2.2]{Friedrich-Schmidt:15}).  Let $\Omega_\rho = \lbrace x\in\Omega: \dist(x, \partial \Omega) > \BBB \rho \EEE \rbrace$ for $\rho>0$. Recall \eqref{rig-eq: SBVfirstdef}, \eqref{rig-eq: Griffith en} and introduce \BBB an auxiliary \EEE energy functional by
\begin{align}\label{rig-eq: Griffith en2}
E_\eps^\rho(y,U) =  \frac{1}{\eps}\int_U W(\nabla y(x)) \,dx + \int_{J_y \cap U} f_\eps^\rho(|[y](x )|)\,d{\cal H}^1(x).
\end{align}
for  $\rho > 0$, $\eps > 0$ and $U \subset \Omega$, where $f_\eps^\rho(\BBB t \EEE ) := \min\lbrace\frac{\BBB t \EEE}{\sqrt{\eps}\rho} ,1 \rbrace$.  Recall the definition $e(G) = \frac{G + G^T}{2}$ for all $G \in \R^{2 \times 2}$.

\begin{theorem}\label{rig-th: rigidity}
Let $\Omega \subset \R^2$  open, bounded with Lipschitz boundary. Let $M>0$ and $0 < \eta, \rho < 1$. Then there are a universal constant $c>0$, constants $\bar{C}=\bar{C}(\Omega,M,\eta)>0$, \BBB $\hat{C}=\hat{C}(\Omega,M,\eta,\rho)>0$, \EEE  and \BBB $\eps_0 = \eps_0(M,\eta,\rho)>0$ \EEE such that the following holds  for all \BBB $0 < \eps \le \eps_0$: \EEE \\
For each $y \in SBV_{M}(\Omega)$ with ${\cal H}^{1}(J_y) \le M$ and $\int_\Omega \dist^2(\nabla y,SO(2) )  \le M\eps$ there is an open set $\Omega_y \subset \Omega$  and a modification  $\hat{y} \in SBV_{cM}(\Omega)$  satisfying 
\begin{align}\label{rig-eq: energy le}
\begin{split}
\BBB (i) \EEE & \ \ \Vert \hat{y} - y \Vert^2_{L^2(\Omega_y)} +  \Vert \nabla \hat{y} - \nabla y \Vert^2_{L^2(\Omega_y)}\le \bar{C}\eps\rho, \ \ \ \ |\Omega\setminus\Omega_y| \le \bar{C}\rho,\\
(ii) & \ \ E_\eps^\rho(\hat{y},\Omega_\rho) \le E_\eps(y) + \bar{C}\rho
\end{split}
\end{align}
with the following properties: We find  a Caccioppoli partition ${\cal P} = (P_j)_j$ of $\Omega_\rho$ with $\sum_j \mathcal{H}^1(\partial^* P_j\cap \Omega_{\rho}) \le \bar{C}$ and for each $P_j$ a corresponding rigid motion $R_j  \cdot +b_j$, $R_j \in SO(2)$ and $b_j \in \R^2$, such that the function $u: \Omega \to \R^2$ defined by
$$
u(x) := \begin{cases} \hat{y}(x) - (R_j\,x +b_j) & \ \ \text{ for } x \in P_j \\
                      0                      & \ \ \text{ for } x \in \Omega \setminus \Omega_\rho \end{cases}
$$
satisfies the estimates
\begin{align}\label{rig-eq: main properties2}
\begin{split}
(i) & \ \, {\cal H}^{1}(J_u) \le \bar{C}, \ \  \ \ \  \ \ \  \ \ \ \   \ \ \   \ \  \ \ \  \ \ \  \ \  \ \ (ii) \,  \ \Vert u\Vert^2_{L^2(\Omega_\rho)} \le \hat{C}\eps, \\
(iii) & \ \, \sum\nolimits_j \Vert e(R^T_j \nabla u)\Vert^2_{L^2(P_j)} \le \hat{C}\eps,  \  \ \ \  \  \ \ \,  (iv)  \ \, \Vert \nabla u\Vert^2_{L^2(\Omega_\rho)} \le \hat{C}\eps^{1-\eta}. 
\end{split}
\end{align}

\end{theorem}

We remark that we get a sufficiently strong bound only for the symmetric part of the gradient (see \eqref{rig-eq: energy le}(iii)) which is not surprising due to the fact that there is no direct analogue of Korn's inequality in $SBV$. However, there is at least a weaker bound on the full absolutely continuous part of the gradient $\nabla u$ (see \eqref{rig-eq: energy le}(iv)) which will essentially be needed to estimate the elastic part of the energy in the passage to the linearized theory (see  \eqref{rig-eq: ff} and \eqref{rig-eq: ff2} below). \BBB In particular,   it will allow us to obtain \eqref{rig-eq: comp2}(ii) after modification of $u$ on a set of  small measure. \EEE

Furthermore, let as briefly note that the uniform bound on the gradient (see \eqref{rig-eq: SBVfirstdef}) in the setting of the nonlinear model is only needed for the application of the rigidity estimate. The condition essentially ensures that the elastic energy cannot concentrate on scales being much smaller than $\eps$. In particular, this is a natural assumption in the investigation of discrete systems, where $\eps$ may be interpreted as the typical interatomic distance.

\begin{rem}\label{rem: rig}
{\normalfont
(i) Estimate \eqref{rig-eq: energy le}(ii) can be refined. Indeed, we obtain 
\begin{align}\label{rig-eq: part + crack}
\begin{split}
(i) & \ \ \BBB \frac{1}{\eps}\int_{\Omega_\rho}W(\nabla \hat{y}(x)) \,dx \le \frac{1}{\eps}\int_{\Omega}W(\nabla y(x)) \, dx + \bar{C}\rho, \EEE \\
(ii) &  \ \ \sum\nolimits_j \frac{1}{2} \mathcal{H}^1(\partial^* P_j \cap \Omega_\rho) + \int_{J_{\hat{y}} \setminus \bigcup_j \partial^* P_j} f_\eps^\rho(|[\hat{y}]|) \,d{\cal H}^1 \le {\cal H}^1(J_y) + \bar{C}\rho.
\end{split}
\end{align}
\BBB We remark that it is indispensable to allow for a  small modification of the deformation in Theorem \ref{rig-th: rigidity} in order to guarantee the sharp energy estimate \eqref{rig-eq: part + crack}(ii). \EEE 

(ii) To derive \eqref{rig-eq: main properties2}\BBB (iii) \EEE one essentially shows
$$
\Vert \nabla u\Vert^4_{L^4(\Omega_\rho)}  = \sum\nolimits_j\Vert  \nabla \hat{y}- R_j\Vert^4_{L^4(P_j)} \le \hat{C}\eps.
$$
The claim then follows from \BBB \eqref{eq:W}, \eqref{rig-eq: part + crack}(i), \EEE $\Vert\dist( \BBB \nabla y, \EEE SO(2))\Vert^2_{L^2(\Omega)} \le M\eps$ and the linearization formula \BBB (see \cite[(3.20)]{FrieseckeJamesMueller:02}) \EEE
\begin{align}\label{rig-eq: linearization}
|e(R^T G -\Id)| =  \dist(G,SO(2)) + \BBB O\EEE(|G- R|^2)
\end{align}
 for $G \in \R^{2 \times 2}$ and $R \in SO(2)$, where $\Id$ denotes the identity matrix.

}
\end{rem}

%----------------------------------------------------------------------------------------------------------------------
\section{Compactness of rescaled configurations}\label{rig-sec: sub, comp1}
%-----------------------------------------------------------------------------------------------------------------------

This section is devoted to the proof of the main compactness result given in Theorem \ref{rig-th: comp1}. Moreover, we also show that  Theorem \ref{rig-th: comp1} provides an alternative proof of the piecewise rigidity result stated in Theorem \ref{rig-cor: cgp}.

\subsection{Preparations} 

For the compactness theorem in $GSBD$  (see Theorem \ref{rig-th: GSBD comp}) it is necessary that the integral for some integrand $\psi$ with $\lim_{\BBB t \EEE \to \infty} \psi( \BBB t \EEE ) = \infty$ is uniformly bounded. We first give a simple criterion for the existence of such a function which is, loosely speaking, based on the condition that the functions coincide in a certain sense on the bulk part of the domain.

\begin{lemma}\label{rig-lemma: concave function}
For every increasing sequence $(b_i)_i \subset (0,\infty)$ with $b_i \to \infty$  there is an increasing concave function $\psi:[0,\infty) \to [0,\infty)$ with $\lim_{t \to \infty} \psi(t) = \infty$ and $\psi(b_i) \le 2^{i}$ for all $i \in \N$.
\end{lemma}

\Proof  Let $f:[0,\infty) \to [0,\infty)$ be the function with $f(0) = 0$, $f(b_i) = 2^{i}$ which is affine on each segment $[b_i,b_{i+1}]$. Clearly, $f$ is increasing and satisfies $f(t) \to \infty$ for $t \to \infty$, but is possibly  not concave. We now construct $\psi$ and first let  $\psi = f$ on $[0,b_1]$. Assume $\psi$ has been defined on $[0,b_i]$ and that \BBB $\psi$ is increasing, concave, satisfies $\psi \le f$ and \EEE $\psi(b_i) = f(b_i) = 2^{i}$. If $f'(b_i-)\ge f'(b_i+)$, we set $\psi = f$ on $[b_{i},b_{i+1}]$. Here, $f'(t\pm)$ denote the one-sided limits of the derivative at point $t$. \BBB This implies that $\psi$ is concave on $[0,b_{i+1}]$ since $\psi'(b_i -) \ge f'(b_i-)$. \EEE

 Otherwise, we let $\psi(t) = f(b_i) + f'(b_i-)(t-b_i)$ for $t \in [b_i,\bar{t}]$, where $\bar{t}$ is the smallest value larger than $b_i$ such that $f(\bar{t}) = f(b_i) + f'(b_i-)(\bar{t}-b_i)$. If $\bar{t}$ does not exist, we are done. If $\bar{t}$ exists, we assume $\bar{t} \in (b_{j-1},b_j]$ and define $\psi=f$ on $[\bar{t}, b_j]$. \BBB Note that if $\bar{t} \in (b_{j-1},b_j)$, we have  $\psi'(\bar{t}-) \ge \psi'(\bar{t}+)$ since $f$ is affine on $[b_{j-1},b_j]$ and $\psi < f$ on $(b_i,\bar{t})$. Thus, $\psi$ is concave on $[0,b_j]$.
 
 Repeating the construction  \EEE we end up with an increasing concave function $\psi$ with $\psi \le f$ and $\psi(t) \to \infty$ for $t \to \infty$. \eop  

\begin{lemma}\label{rig-lemma: concave function2}
Let $\Omega \subset \R^2$ and let $(y^l)_l \subset L^1(\Omega)$   be a sequence satisfying $|\Omega \setminus \bigcup_{n\in \N} \bigcap_{l \ge n} \lbrace |y^n - y^l| \le 1 \rbrace|=0$. Then there is a not relabeled subsequence such that 
$$\int_{\Omega}\psi(|y^l|)\le C$$
for a constant $C>0$ independent of $l$, where $\psi$ is an increasing continuous function  with $\lim_{t \to \infty} \psi(t) =  \infty$.
\end{lemma}

\Proof  Define $C_l := \max_{1 \le i \le l} \Vert y^i \Vert_{L^1(\Omega)}$ for all $l \in \N$. Let $A_n = \bigcap_{l \ge n} \lbrace |y^n - y^l|\le 1 \rbrace$ and set $B_1 = A_1$ as well as $B_n = A_n \setminus \bigcup^{n-1}_{m=1} B_{m}$ for all $n \in \N$. The sets $(B_n)_n$ are pairwise disjoint with $\sum_n |B_n| = |\Omega|$. We choose $0=n_1 < n_2 < \ldots$ such that $\sum_{1\le n \le n_{i}} \frac{|B_n|}{|\Omega|} \ge 1 - 4^{-i}$. We let $B^i = \bigcup^{n_{i+1}}_{n=n_i+1} B_n$ and observe $|B^i|\le  4^{-i}|\Omega| $. 

We pass to the  subsequence  $(n_i)_i \subset \N$ and choose $ E^i \supset B^i$ such that $| E^i| = 4^{-i} |\Omega|$. Let $b_i = \frac{C_{n_{i+1}}}{| E^i|} + 2=  4^{i}\frac{C_{n_{i+1}}}{|\Omega|}+ 2$ for $i\in \N$  and note that $(b_i)_i$ is increasing with $b_i \to \infty$.  By Lemma \ref{rig-lemma: concave function} we get an increasing concave function $\psi:[0,\infty) \to [0,\infty)$ with $\lim_{t \to \infty} \psi(t) = \infty$ and $\psi(b_i) \le 2^{i}$ for all $i \in \N$. Clearly, $\psi$ is also continuous.   

For $\hat{B}^{ i} := \Omega \setminus \bigcup^{n_i}_{n=1} B_n$ we have $|\hat{B}^{ i}| \le 4^{-i}|\Omega|$ and choose $\hat{E}^{ i} \supset \hat{B}^i$ with $|\hat{E}^{ i}| = 4^{-i} |\Omega|$. We then obtain  $\frac{C_{n_i}}{|\hat{E}^{ i}|}=  4^{i}\frac{C_{n_{i}}}{|\Omega|} \le b_{i}$.  Now let $l = n_i$. Using Jensen's inequality, the definition of the sets $B^i$, $\Vert y^l \Vert_{L^1(\Omega)} \le C_l$  and the monotonicity of $\psi$ we compute
\begin{align}\label{rig-eq: psi est}
\begin{split}
\int_\Omega \psi(|y^l|) &=\sum\nolimits_{1 \le j \le i-1} \int_{B^j} \psi(|y^l|) + \int_{\hat{B}^{ i}} \psi(|y^l|)   \\
& \BBB \le \EEE \sum\nolimits_{1 \le j \le i-1}  \int_{B^j} \psi(|y^{n_{j+1}}|+ 2) + \int_{\hat{B}^{i}} \psi(|y^l|) \\
& \le \sum\nolimits_{1 \le j \le i-1}  |E^j| \psi\Big(\dashint_{E^j} |y^{n_{j+1}}|+2\Big) + |\hat{E}^{ i}| \psi\Big(\dashint_{\hat{E}^{ i}} |y^l|\Big) \\
& \le \sum\nolimits_{1 \le j \le i-1} 4^{-j}|\Omega| 2^{j} + 4^{-i}|\Omega| 2^{i} \le |\Omega|\sum\nolimits_{j \in \N} 2^{-j}.
\end{split}
\end{align}   
As the estimate is independent of $l \in (n_i)_i$, this yields $\int_\Omega \psi(|y^l|) \le C$ uniformly in $l$, as desired. \eop

\subsection{Proof of Theorem \ref{rig-th: comp1}}

Now we are in a position to give the proof of the main compactness result. In the first part we show that \eqref{rig-eq: comp2}-\eqref{rig-eq: comp1-2} hold.

\noindent {\em Proof of Theorem \ref{rig-th: comp1}, part 1.}  Let $(\eps_k)_k$ be a sequence \BBB with $\eps_k \to 0$. \EEE Let $y_k \in SBV_M(\Omega)$ with $E_{\eps_k}(y_k) \le C$ be given. \BBB Possibly passing to a larger $M$, \EEE we get $\Vert \dist (\nabla y_k,SO(2)) \Vert^2_{L^2(\Omega)} \le  M\eps_k$ \BBB by \eqref{eq:W} and  ${\cal H}^1(J_{y_k}) \le M$ for all $k \in \N$. \BBB In the following generic constants only  depending on $\Omega$ and $M$ will be denoted by $C$. \EEE

\smallskip
\BBB\emph{Step I:} \EEE Choose $\rho_0 >0$ small and let $\rho_l = 2^{- 3l}  \rho_0$ for all $l\in \N$. \BBB  We apply Theorem \ref{rig-th: rigidity} for $\rho = \rho_l$ and $\eta=\frac{1}{5}$ (the choice of $\eta$ is related to the exponent $-\frac{1}{8}$ in \eqref{rig-eq: comp2}(ii)). Denote by $c$, $\bar{C}=\bar{C}(\Omega,M,\eta)$,  $\hat{C}_l=\hat{C}_l(\Omega,M,\eta,\rho_l)$ the constants in Theorem \ref{rig-th: rigidity}.  For each $l \in \N$ there exists $\kappa_l = \kappa_l(M,\eta,\rho_l)$ such that for $k \ge \kappa_l$ \EEE we find  modifications  $y^l_k \in SBV_{cM}(\Omega, \R^2)$ with $E_{\eps_k}^{\rho_l} (y^l_k,  \Omega_{\rho_l}) \le E_{\eps_k}(y_k) + \bar{C}\rho_l$ and
\begin{align}\label{rig-eq: approx en}
 \Vert y^l_k -  y_k \Vert^2_{L^2(\Omega^l_k)} + \Vert \nabla y^l_k - \nabla  y_k \Vert^2_{L^2(\Omega^l_k)} \le \bar{C}\eps_k\rho_l,
\end{align}
where $\Omega^l_k := \Omega_{y^l_k}$ with $|\Omega \setminus \Omega^l_k| \le \bar{C}\rho_l$. We further get Caccioppoli partitions $(P^{k,l}_j)_j$ of $\Omega_{\rho_l}$ with $\sum_j \mathcal{H}^1(\partial^* P^{k,l}_j \cap \Omega_{\rho_l}) \le \bar{C}$ and corresponding piecewise rigid motions $T_k^l:= \sum_j(R^{k,l}_j \cdot + b^{k,l}_j) \chi_{P^{k,l}_j} \BBB + \id\chi_{\Omega \setminus \Omega_{\rho_l}} \EEE $ such that the functions $v^l_k : \Omega \to \R^2$ defined by
\begin{align}\label{rig-eq: comp13}
v^l_k (x)  = \begin{cases}
\frac{ 1}{\sqrt{\eps_k}} (R^{k,l}_j)^T \big( y^l_k(x) - (R^{k,l}_j\,x + b^{k,l}_j)\big) & \text{ for } x \in P^{k,l}_j, \ j \in \N, \\
0 & \text{ else,}
\end{cases}
\end{align}
satisfy  by \eqref{rig-eq: main properties2}
\begin{align}\label{rig-eq: comp11}
{\cal H}^1(J_{v^l_k}) \le \bar{C}, \ \ \  \Vert v^l_k \Vert_{L^2(\Omega)} + \Vert e (\nabla v^l_k) \Vert_{L^2(\Omega)} \le \hat{C}_l, \ \ \ \Vert \nabla v^l_k \Vert^2_{L^2(\Omega)} \le \hat{C}_l\eps_k^{-\BBB 1/5 \EEE }.
\end{align}
\BBB We recall $\Vert y^l_k \Vert_\infty \le cM$ for all $k \ge \kappa_l$. Thus, possibly passing to other (not relabeled) constants $b^{k,l}_j$ in \eqref{rig-eq: comp13}, we can assume that $|b^{k,l}_j| \le CM$ for $C=C(\Omega,c)$ and that \eqref{rig-eq: comp11} still holds. \EEE Each partition may be extended to $\Omega$ by adding the element $\Omega \setminus \Omega_{\rho_l}$. \BBB  As for $\rho_0$  small  enough we get  ${\cal H}^1(\partial \Omega_{\rho_l}) \le C{\cal H}^1(\partial \Omega)$ (see \cite[Theorem 4.1]{Doktor}) for all $l \in \N$, \EEE there is   $C=C(\bar{C},\Omega)$ such that  
\begin{align}\label{eq: new last}
\sum\nolimits_j \mathcal{H}^1(\partial^* P^{k,l}_j) \le C.
\end{align}

\BBB\emph{Step II:} \EEE Using a diagonal argument we get a (not relabeled) subsequence of $(\eps_k)_k$ such that by Theorem \ref{rig-th: GSBD comp} for every $l \in \N$ we find   $v^l \in GSBD^2(\Omega)$  with 
\begin{align}\label{rig-eq: T conv3}
v^l_k \to v^l \text{ a.e. in } \Omega \text{  \ \ \   and  \ \ \    } e(\nabla v^l_k) \rightharpoonup e(\nabla v^l) \text{ weakly in } L^2(\Omega, \R^{2 \times 2}_{\rm sym})
\end{align}
for $k \to \infty$. By \BBB Theorem \ref{th: comp cacciop}, Theorem \ref{th: piecewise const}, \eqref{eq: new last}  and the fact that  $|b^{k,l}_j| \le CM$ \EEE we obtain an  (ordered) partition $(P^l_j)_j$ of $\Omega$ with $\sum_j \mathcal{H}^1(\partial^* P^l_j) \le C$ and a piecewise rigid motion $T^l := \sum_j(R^l_j \cdot + b^l_j)\chi_{P^l_j}$  such that  for all $l \in \N$ we get (again up to a subsequence) \BBB $|P^{k,l}_j\triangle P^l_j| \to 0$, $R^{k,l}_j \to R^l_j $, and  $b^{k,l}_j \to b^l_j $ \EEE for all $j \in \N$ as   $k \to \infty$. This also implies 
\begin{align}\label{rig-eq: T conv2}
 \sum\nolimits_j |P^{k,l}_j \triangle P^{l}_j| +  \Vert T_k^l - T^l \Vert_{L^2(\Omega)} + \Vert \nabla T_k^l - \nabla T^l \Vert_{L^2(\Omega)}\to 0
\end{align}
for $k \to \infty$.  We now show that
\begin{align}\label{rig-eq: comp12}
\Vert v^l\Vert_{L^1(\Omega)}  \le C\Vert v^l \Vert_{L^2(\Omega)} \le C \hat{C}_l, \ \ \ {\cal H}^1(J_{v^l}) \le \bar{C}, \ \ \ \Vert e(\nabla v^l)\Vert^2_{L^2(\Omega)} \le C.
\end{align}
The first two claims follow directly from \eqref{rig-eq: comp11} and \eqref{rig-eq: convergence sense}. To see the third estimate, we let $\BBB \phi^l_k (x) \EEE:= \chi_{[0,\eps_k^{-1/8}]} (|\nabla v^l_k (x)|)$ \BBB (cf. \eqref{rig-eq: comp2}(ii)). \EEE Moreover,  we obtain by an elementary computation (cf. \eqref{rig-eq: linearization}) $\dist^2(G,SO(2)) = | e(R^T G - \Id)|^2 + \omega_{\rm dist  }(R^TG-\Id)$ for $G \in \R^{2 \times 2}$, $R \in SO(2)$ with $\sup\lbrace |G|^{-3}\omega_{\rm dist}(G): |G| \le 1\rbrace \le C$. We compute  by \eqref{eq:W} and \eqref{rig-eq: comp13}
\begin{align}\label{rig-eq: ff}
\begin{split}
C & \ge E^{\rho_l}_{\eps_k} (y^l_{k},\Omega_{\rho_l})  \geq \frac{C}{\eps_k} \int_{\Omega_{\rho_l}} \dist^2(\nabla y^l_k,SO(2) )\\ & \ge \frac{C}{\eps_k} \sum\nolimits_j\int_{P^{k,l}_j \BBB \cap \Omega_{\rho_l} \EEE } \phi^l_k \Big(|{e}( (R^{k,l}_j)^T\nabla y^l_k - \Id)|^2 + \omega_{\rm dist}((R^{k,l}_j)^T\nabla y^l_k - \Id)  \Big) \\ & = C \int_{\Omega} \phi^l_k \Big(|e(\nabla v^l_k)|^2 + \frac{1}{\eps_k}\omega_{\rm \dist}(\sqrt{\eps_k} \nabla v^l_k )  \Big). 
\end{split}
\end{align}
The second term of the integral can be estimated by 
\begin{align}\label{rig-eq: ff2}
\int_{\Omega} \phi^l_k \frac{1}{\eps_k}\omega_{\rm dist}(\sqrt{\eps_k} \nabla v^l_k )   = \int_{\Omega}\phi^l_k \sqrt{\eps_k} |\nabla v^l_k |^3 \frac{\omega_{\rm dist}(\sqrt{\eps_k} \nabla v^l_k) }{|\sqrt{\eps_k} \nabla v^l_k |^3} \le C\eps_k^{\frac{1}{8}} \to 0.
\end{align}
As $e(\nabla v^l_k) \rightharpoonup e(\nabla v^l)$ weakly in $L^2(\Omega)$ and $\phi^l_{k} \rightarrow 1$ boundedly in measure on $\Omega$  by \eqref{rig-eq: comp11}, it follows $\phi^l_{k} e(\nabla v^l_k) \rightharpoonup e(\nabla v^l)$ weakly in $L^2(\Omega)$. By lower semicontinuity we obtain $\Vert e(\nabla v^l)\Vert^2_{L^2(\Omega)} \le C$ for a constant particularly independent of $\rho_l$ which concludes \eqref{rig-eq: comp12}.

\smallskip
\BBB\emph{Step III:} \EEE We now want to pass to the limit $l \to \infty$. Similarly as  in the argumentation leading to \eqref{rig-eq: T conv2}, by the compactness result for piecewise constant functions (see Theorem \ref{th: piecewise const}) we find a partition $(P_j)_j$ of $\Omega$ and a piecewise rigid motion $T := \sum_j(R_j \cdot +b_j)\chi_{P_j}$ such that for a suitable (not relabeled) subsequence 
\begin{align}\label{rig-eq: T conv1}
 \sum\nolimits_j |P^{l}_j \triangle P_j| + \Vert T^l - T \Vert_{L^2(\Omega)} + \Vert \nabla T^l - \nabla T \Vert_{L^2(\Omega)}\to 0
\end{align}
for $l \to \infty$. Recalling \eqref{rig-eq: T conv2}  and using a diagonal argument we can choose a (not relabeled) subsequence of $(\rho_l)_l$ and afterwards of $(\eps_k)_k$ such that  for all $l$ we have
\begin{align}\label{rig-eq: subsubseq}
\sum\nolimits_j |P^l_j \triangle P_j| \le 2^{-l}, \ \ \ \sum\nolimits_j |P^{k,l}_j \triangle P^{l}_j| \le 2^{-l} \ \text{ for all } \ k \ge l.
\end{align} 
We see that  the  compactness result in $GSBD$ cannot be  applied directly on the sequence $(v^l)_l$ as the $L^2$ bound in \eqref{rig-eq: comp12} depends on $\rho_l$. We now show that by choosing the rigid motions on the elements of the partitions appropriately (see \eqref{rig-eq: comp13}) we can construct the sequence $(v^l)_l$ such that we obtain 
\begin{align}\label{rig-eq: Asetprep}
\Big|\Omega \setminus \bigcup\nolimits_{n \in \N} \bigcap\nolimits_{ m \ge n}\lbrace |v^n - v^{m}| \le 1\rbrace\Big|=0
\end{align}
and thus Lemma  \ref{rig-lemma: concave function2} is applicable. 

\smallskip
\BBB\emph{Step IV:} \EEE  We fix $k \in \N$ and describe an iterative procedure to redefine $R^{k,l}_j, b^{k,l}_j$ for all $l$ with \BBB $k \ge \kappa_l$ \EEE and $j \in \N$. Let $\tilde{v}^1_k= {v}^1_k$ as defined in \eqref{rig-eq: comp13} and assume \BBB $\tilde{v}^l_k$ with corresponding \EEE $\tilde{R}^{k,l}_j,  \tilde{b}^{k,l}_j$ have been chosen  (which possibly differ from $R^{k,l}_j,  b^{k,l}_j$)  such that \eqref{rig-eq: comp11} still holds possibly passing to a larger constant \BBB $\tilde{C}_l= \tilde{C}_l(\Omega,M,\eta,l)$. Let $I^{k,l}_1 = \lbrace j: |P^{k,l+1}_{j} \cap P^{k,l}_{j}| \ge 4\bar{C}\rho_l\rbrace$ and $I^{k,l}_2 = \N \setminus I^{k,l}_1$. \EEE Define
\begin{align}\label{new3}
\begin{split}
&\tilde{R}^{k,l+1}_{j} = \tilde{R}^{k,l}_{j}, \ \ \ \ \ \  \tilde{b}^{k,l+1}_{j} = \tilde{b}^{k,l}_{j} \ \ \ \ \ \,  \text{  for } j \in I^{k,l}_1,\\
&\tilde{R}^{k,l+1}_{j} = R^{k,l+1}_{j}, \ \ \ \tilde{b}^{k,l+1}_{j} = b^{k,l+1}_{j} \ \ \ \text{  for } j \in I^{k,l}_2.
\end{split}
\end{align}
\BBB Consider $j \in I^{k,l}_1$ and define $R_j' = {R}^{k,l+1}_{j} - \tilde{R}^{k,l}_{j}$, $b_j' = {b}^{k,l+1}_{j} - \tilde{b}^{k,l}_{j}$, $P'_j = P^{k,l+1}_{j} \cap P^{k,l}_{j} \cap \Omega^l_k \cap \Omega_k^{l+1}$ for shorthand. By the triangle inequality,  \eqref{rig-eq: approx en} and \eqref{rig-eq: comp13}  we get
\begin{align}\label{eq:NEW}
\begin{split}
\Vert R_j'  \cdot + b_j' \Vert_{L^2(P'_j)} &\le \sqrt{\eps_k}(\Vert \tilde{v}^l_k \Vert_{L^2(\Omega)} + \Vert {v}^{l+1}_k \Vert_{L^2(\Omega)}) + \Vert y^l_k - y^{l+1}_k \Vert_{L^2(\Omega^l_k \cap \Omega_k^{l+1})}  \\
& \le  \sqrt{\eps_k} (\tilde{C}_l + \hat{C}_{l+1} + C\sqrt{\rho_l}+ C\sqrt{\rho_{l+1}}) \le C'_l \sqrt{\eps_k}
\end{split}
\end{align}
for a constant $C'_l = C'_l(\Omega,M,\eta, l)$, where in the penultimate step we used that \eqref{rig-eq: comp11} holds for $\tilde{v}^l_k$ and ${v}^{l+1}_k$. Herefrom we now derive $|R_j'| \le C'_l \sqrt{\eps_k}$. Indeed, if $R_j' \neq 0$, then $R_j'$ is invertible and a short computation yields
\begin{align}\label{eq:NEW2}
\tfrac{1}{\sqrt{2}}|R_j'|\Vert \cdot - z\Vert_{L^2(P'_j \setminus B_\lambda(z))} \le \Vert R_j'  \cdot + b_j' \Vert_{L^2(P'_j)}  \le C'_l \sqrt{\eps_k},
\end{align}
 where $z:= - (R_j')^{-1}b_j'$ and  $B_\lambda(z)$ denotes the ball with center $z$ and radius $\lambda = (\pi^{-1}\bar{C}\rho_l)^{1/2}$. Then by definition of $I^{k,l}_1$ and $|\Omega \setminus \Omega_k^l| \le \bar{C}\rho_l$ we find $|P_j' \setminus B_\lambda(z)| \ge |P^{k,l+1}_{j} \cap P^{k,l}_{j}| - |\Omega \setminus \Omega_k^l| - |\Omega \setminus \Omega_k^{l+1}| - |B_\lambda(z)| \ge \bar{C}\rho_l$, which together with \eqref{eq:NEW2} implies the claim for $C'_l$ sufficiently large.

Recalling \eqref{eq:NEW} we then also find $|b_j'| \le C'_l \sqrt{\eps_k}$ and summing over all components we derive
$$\sum_{j \in I^{k,l}_1} \Big( \Vert R_j' \Vert^2_{L^2(P_j^{k,l+1})} + \Vert R_j' \Vert^4_{L^4(P_j^{k,l+1})} + \Vert R_j' \, \cdot +  b_j'\Vert^2_{L^2(P_j^{k,l+1})} \Big) \le \# I^{k,l}_1 C'_l \eps_k \le  \frac{|\Omega|C'_l}{4\bar{C}\rho_l}  \eps_k,$$
where in the last step we used the definition of $I^{k,l}_1$. Define $\tilde{v}_k^{l+1}$ as in \eqref{rig-eq: comp13} with $\tilde{R}^{k,l+1}_j,  \tilde{b}^{k,l+1}_j$ instead of ${R}^{k,l+1}_j, {b}^{k,l+1}_j$. The previous estimate together with the fact that \eqref{rig-eq: comp11} holds for ${v}^{l+1}_k$ now shows \eqref{rig-eq: comp11}  for $\tilde{v}_k^{l+1}$.  Indeed, the estimates for $\Vert \tilde{v}^{l+1}_k \Vert_{L^2(\Omega)}$, $\Vert \nabla \tilde{v}^{l+1}_k \Vert_{L^2(\Omega)}$ follow directly and for $\Vert e (\nabla \tilde{v}^{l+1}_k) \Vert_{L^2(\Omega)}$ we argue as in Remark \ref{rem: rig}(ii). \EEE

\BBB Note that  as $(\hat{C}_l)_l$ also   $(\tilde{C}_l)_l$ converges to infinity. For simplicity the modified functions and rigid motions will still be denoted by $v^l_k$, $R_j^{k,l}$ and $b_j^{k,l}$ in the following.  
By a diagonal argument we can choose  a further (not relabeled) subsequence of $(\eps_k)_k$ such that the modifications $v^l_k$ exist for all $l \in \N$ and $k \ge l$\EEE.

\smallskip
\BBB\emph{Step V:} \EEE We define $ \BBB A^n_{k,l} \EEE = \bigcap_{n \le m \le l} \lbrace |v^m_k - v^n_k| \le \frac{1}{2} \rbrace$ for all $n \in \N$ and  $n \le l \le k$. If we show
\begin{align}\label{rig-eq: Akl2}
|\Omega \setminus A^n_{k,l}| \le C2^{-n},
\end{align}
then \eqref{rig-eq: Asetprep} follows. Indeed, for given $l\ge n$ we can choose $K=K(l)\ge l$ so large that $|\lbrace |v^m_K - v^m| > \frac{1}{4} \rbrace| \le 2^{-m}$ for all $n \le m \le l$ since $v^m_k \to v^m$ in measure for $k \to \infty$. This implies 
$$\big|\Omega \setminus \bigcap\nolimits_{n \le m \le l} \lbrace |v^m -v^n|\le 1 \rbrace\big|\le |\Omega \setminus A^n_{K,l}| + \sum\nolimits_{n \le m \le l} |\lbrace |v^m_K - v^m| > \text{\scriptsize $\frac{1}{4}$} \rbrace| \le C2^{-n}.$$
Passing to the limit $l \to \infty$ we find $|\Omega \setminus \bigcap\nolimits_{ m \ge n} \lbrace |v^{m} -v^n| \le 1 \rbrace| \le C2^{-n}$ and taking the union over all $n \in \N$ we derive \eqref{rig-eq: Asetprep}.

 To show \eqref{rig-eq: Akl2} we proceed in two steps.  Employing the redefinition of the piecewise rigid motions we first show that the set where $T^{m}_k, n \le m \le l$, differ is small. Afterwards, we use \eqref{rig-eq: approx en} to find that the set where $y^{m}_k, n \le m \le l$, differ is small. We define $ \BBB B^n_{k,l} \EEE = \bigcap_{n \le m \le l} \lbrace T^m_k = T^n_k \rbrace$ for $k \ge l  \ge n$ and prove that 
\begin{align}\label{rig-eq: Bkl}
|\Omega \setminus B^n_{k,l}| \le C2^{-n}
\end{align}
for all $k \ge l  \ge n$.  To this end, consider $\lbrace T^m_k = T^{m+1}_k \rbrace$ for $n \le m \le l-1$ and first note that by \eqref{rig-eq: subsubseq} we have $\sum_j |P_j^{k, m+1} \triangle P_j^{k,m}| \le 3 \cdot 2^{-m}$. Define  $J_1 \subset \N$ such that \BBB $|P_j^{k,m+1}| \le 8\bar{C}\rho_m$ \EEE for all $j \in J_1$ and let $J_2 \subset \N \setminus J_1$ such that $|P^{k, m+1}_j \cap P^{k, m}_j| > \frac{1}{2} |P_j^{k, m+1}|$ for all  $j \in J_2$. Observe that $|P_j^{k, m+1}| \le 2 |P_j^{k, m+1} \setminus P_j^{k,m}|$ for  $j \in J_3 := \N \setminus (J_1 \cup J_2)$. \BBB Using the isoperimetric inequality and the the fact that $(\rho_m)_m \subset (2^{- 3m}\rho_0)_m$  we find by \eqref{eq: new last}
\begin{align*}
\sum_{j \in J_1} |P_j^{k,m+1}| \le (8\bar{C}\rho_m)^{\frac{1}{2}} \sum_{j \in J_1} |P_j^{k,m+1}|^{\frac{1}{2}} \le C2^{-m}\sum\nolimits_j \mathcal{H}^1(\partial^* P_j^{k,m+1})\le C2^{-m}.
\end{align*}
\EEE Due to the above construction of the rigid motions \BBB (see \eqref{new3}) \EEE we obtain $\lbrace T^m_k = T^{m+1}_k \rbrace \supset \bigcup_{j \in  J_2} (P^{k,m+1}_j \cap P^{k,m}_j)$  and therefore  
\begin{align*}
|\Omega \setminus \lbrace T^m_k = T^{m+1}_k \rbrace| & \le \sum\nolimits_{j \in J_2} |P_j^{k, m+1} \setminus P_j^{k, m}| + \sum\nolimits_{j \in J_1 \cup J_3} |P_j^{k, m+1}| \\
&\le \sum\nolimits_{j \in J_2} |P_j^{k, m+1} \setminus P_j^{k, m}| + \sum\nolimits_{j \in J_3} 2|P_j^{k, m+1} \setminus P_j^{k, m}| \BBB + C2^{-m} \EEE \\
& \le \BBB 2\sum\nolimits_j |P_j^{k, m+1} \triangle P_j^{k,m}| \EEE +  C2^{-m} \le C2^{-m}.
\end{align*}
Summing over $n \le m \le l-1$ we establish \eqref{rig-eq: Bkl}.   Now recalling \eqref{rig-eq: approx en}, \eqref{rig-eq: Bkl},  $|\Omega \setminus \Omega_k^l| \le \bar{C}\rho_{l}$ and the fact that $(\rho_l)_l \subset (2^{- 3l}\rho_0)_l$   we find 
$$|\Omega \setminus A^n_{k,l}| \le  |\Omega \setminus B^n_{k,l}| + \sum\nolimits_{n \le m \le l-1} |\lbrace |y^{m+1}_k - y^m_k| > 2^{-m  -1} \sqrt{\eps_k} \rbrace| \le C2^{-n}$$ 
for all $k \ge l \ge n$,  as desired. 

\smallskip
\BBB\emph{Step VI:} \EEE By \eqref{rig-eq: comp12}  and \eqref{rig-eq: Asetprep} we can apply Lemma \ref{rig-lemma: concave function2} on the sequence $(v^l)_l$. We employ Theorem \ref{rig-th: GSBD comp} and obtain a function $v \in GSBD(\Omega)$ and a further not relabeled subsequence with $v^l \to v$  a.e in $\Omega$ and $e(\nabla v^l) \rightharpoonup e(\nabla v)$ weakly in $L^2(\Omega, \R^{2 \times 2}_{\rm sym})$.

We now select a suitable diagonal sequence such that \BBB \eqref{rig-eq: comp2}-\eqref{rig-eq: comp1-2} hold. \EEE First, we may suppose that after an infinitesimal modification we have $v_k^l \in W^{2,\infty}(\Omega \setminus \overline{J_{v^l_k}})$ (see \cite{Cortesani-Toader:1999}). Consequently, by the coarea formula \cite[ Theorem 3.40]{Ambrosio-Fusco-Pallara:2000} we get $\mathcal{H}^1(\partial^*\lbrace |\nabla v^l_k| \le \lambda\eps_k^{-1/8} \rbrace) < \infty$ for all $\lambda \in (\frac{1}{2},1)\setminus H_k^l$, where $H_k^l$ is an $\mathcal{L}^1$-negligible set. Choosing $\lambda \in (\frac{1}{2},1)\setminus \bigcup_{k,l \in \N}H_k^l$ and defining $\hat{\phi}^l_k = \chi_{[0,\lambda\eps_k^{-1/8}]} (|\nabla v^l_k (x)|)$, the functions  $\hat{v}_k^l :=  \hat{\phi}_k^l  v^l_k$ lie in $SBV(\Omega)$ by  \cite[ Theorem 3.84]{Ambrosio-Fusco-Pallara:2000}. Recalling the definition of $\phi^l_k$ before \eqref{rig-eq: ff}, we observe that by \eqref{rig-eq: ff}, \eqref{rig-eq: ff2} the functions fulfill  $\Vert e(\nabla\hat{v}_k^l)\Vert_{L^2(\Omega)}  \le C$ 
and $\Vert \nabla \hat{v}_k^l \Vert_\infty \le \eps_k^{-1/8}$ for a constant independent of $k,l \in \N$. Moreover, by \eqref{rig-eq: comp11} we get $\hat{\phi}^l_k \to 1$ in measure on $\Omega$ as $k \to \infty$. 

 \EEE As weak convergence in $L^2$ is metrizable on bounded sets and convergence in measure is metrizable (take $(f,g) \mapsto \int_\Omega \min\lbrace |f-g|,1\rbrace$), we can apply a diagonal sequence argument and find a not relabeled subsequence $(y_n)_{n}$ and a corresponding diagonal sequence $(w_n)_{n \in \N} \subset   (\hat{v}_k^l)_{k,l}  $ with corresponding partitions $(P^n_j)_j$ and piecewise rigid motions $(T_n)_n$  such that by \eqref{rig-eq: T conv3}, \eqref{rig-eq: T conv2} and \eqref{rig-eq: T conv1} 
\begin{align}\label{new4}
& w_n \to v \ \text{ in measure on} \ \Omega, \  \ \   e(\nabla w_n) \rightharpoonup e(\nabla v)  \text{ weakly in } L^2(\Omega),  \\
&  T_n \to T \text{ in } L^2(\Omega), \ \ \  \nabla T_n \to \nabla T \text{ in } L^2(\Omega), \  \ \ \BBB |P^n_j \triangle P_j| \to 0 \EEE \ \ \text{for all } j \in \N \notag
\end{align}
for $n \to \infty$. Up to a further subsequence we can assume $w_n \to v$ a.e.  and $\nabla T_n \to \nabla T$ a.e.  in $\Omega$. Finally, define  $u_n = \nabla T_n w_n$  for all $n \in \N$ and let  $u= \nabla T v$. Observe that \eqref{rig-eq: comp2}(ii), \eqref{rig-eq: comp1}, and \eqref{rig-eq: comp1-2} hold. Moreover, as  $ \hat{\phi}_k^l \to 1$  in measure on $\Omega$  and  $|\Omega \setminus \Omega_k^l| \to 0$ for $k,l \to \infty$, we also get \eqref{rig-eq: comp2}(i) recalling  \eqref{rig-eq: approx en}, \eqref{rig-eq: comp13} and possibly passing to a further subsequence.  \eop

To complete the proof of Theorem \ref{rig-th: comp1}, it remains to show \eqref{rig-eq: comp3}. 

\noindent {\em Proof of Theorem \ref{rig-th: comp1}, part 2.}  \BBB To see \eqref{rig-eq: comp3}(i), it suffices to recall that each $\nabla T_n^T u_n$ coincides with some $ \hat{\phi}^l_k v^l_k$ and thus $\Id + \sqrt{\eps_n}\nabla T_n^T \nabla u_n = (\nabla T^l_k)^T\nabla y^l_k$ a.e. on $\lbrace  \hat{\phi}^l_k =1 \rbrace \cap \Omega_{\rho_l}$ by \eqref{rig-eq: comp13}. The assertion then follows from \eqref{rig-eq: part + crack}(i) and the frame indifference of $W$.  We now show \eqref{rig-eq: comp3}(ii). \EEE To this end, the estimate is first carried out in terms of the \BBB auxiliary \EEE functionals (see \eqref{rig-eq: Griffith en2}). Afterwards, we conclude by passing to the limit $\rho \to 0$.

\BBB Let $v$ as given in \eqref{new4} and recall $u = \nabla T v$. \EEE The sets $J_v^{\BBB d \EEE } := \lbrace x \in J_v: [v](x) = d \rbrace$ for $d \in B_1(0)\setminus \lbrace 0 \rbrace$ are pairwise disjoint with ${\cal H}^1$-$\sigma$ finite union, i.e. ${\cal H}^1(J_v^d) = 0$ up to at most countable values of $d$. Consequently, we can choose a sequence $(b_j)_j$ with $0 < |b_j| \le 1$ such that $b_i \neq b_j$ and  ${\cal H}^1(J_v^{b_i - b_j}) = 0$ for $i \neq j$. Replacing $v$ by  $\tilde{v} = v + \sum_j b_j \chi_{P_j}$, we thus obtain ${\cal H}^1 ( (\bigcup_j \partial^* P_j \BBB \cap \Omega ) \EEE \setminus J_{\tilde{v}}) = 0$.  We first show that for \eqref{rig-eq: comp3}(ii) it suffices to prove
\begin{align}\label{rig-eq: compa1}
\liminf\nolimits_{k \to \infty}  {\cal H}^1(J_{y_k} ) \ge {\cal H}^1(J_{\tilde{v}}).
\end{align}
Indeed, we get ${\cal H}^1(J_{\tilde{v}}) = {\cal H}^1(J_{u} \setminus \partial P) + {\cal H}^1(\partial P \cap \Omega)$, \BBB where for shorthand $\partial P = \bigcup_j \partial^* P_j$:  \BBB We have ${\cal H}^1(J_{\tilde{v}}) = {\cal H}^1(J_{\tilde{v}} \cup (\partial P \cap \Omega)) = {\cal H}^1(\partial P \cap \Omega) + {\cal H}^1(J_{\tilde{v}} \setminus \partial P)$. Then it suffices to note  ${\cal H}^1(J_{\tilde{v}} \setminus \partial P) = {\cal H}^1(J_{u} \setminus \partial P)$. \EEE

We now show \eqref{rig-eq: compa1} in two steps first passing to the limit $k \to \infty$ and then letting $l \to \infty$. We replace $v^l_k$ \BBB (see \eqref{rig-eq: comp13}) \EEE by $\tilde{v}^l_k = v^l_k + \sum_j b_j \chi_{P^{ k,l}_j}$ and $v^l$ by $\tilde{v}^l = v^l + \sum_j b_j \chi_{P^{l}_j}$   noting that $\tilde{v}^l_k \to \tilde{v}^l$ for $k \to \infty$ \BBB (cf. \eqref{rig-eq: T conv3}) \EEE and $\tilde{v}^l \to \tilde{v}$ for $l \to \infty$ in the sense of \eqref{rig-eq: convergence sense}.  In the following we write $J^l_k = J_{\tilde{v}^l_k} \cap \Omega_{\rho_l}$ and $\partial  P^{k,l} := \bigcup_j \partial^* P^{k,l}_j$ for shorthand, \BBB where $\Omega_{\rho_l}$ was defined before \eqref{rig-eq: Griffith en2}. \EEE We obtain  by \eqref{rig-eq: comp13}, \eqref{rig-eq: part + crack}(ii) \BBB and Theorem \ref{th: local structure} \EEE
\begin{align}\label{rig-eq: comp1.2}
{\cal H}^1(J_{y_k})  + C\rho_l & \ge \BBB \int_{J_{y^l_k} \setminus \partial P^{k,l}} f^{\rho_l}_{\eps_k}(|[y^l_k]|) \, d{\cal H}^1 + \mathcal{H}^1( \partial  P^{k,l} \cap \Omega_{\rho_l}) \EEE  \\
& \ge \int_{J^l_k\setminus \partial P^{k,l}} \theta_{\rho_l}(|[\tilde{v}^l_k]|) \, d{\cal H}^1 + \mathcal{H}^1( \partial  P^{k,l} \cap \Omega_{\rho_l})   \ge \int_{J^l_k} \theta_{\rho_l}(|[\tilde{v}^l_k]|) \, d{\cal H}^1,\notag
\end{align}
where $\theta_{\sigma}(\BBB t \EEE ) := \min\lbrace \frac{t}{\sigma},1\rbrace$ for $\sigma>0$.  We cannot directly apply  lower semicontinuity results for $GSBD$ functions due to the involved function $\theta_{\rho_l}$. We therefore pass to the limit $k \to \infty$ on one-dimensional sections.

Recall the measure $\hat{\mu}^{\sigma,\xi}_{\tilde{v}^l}$ defined in \eqref{rig-eq: lemma2**} for $\sigma \ge 0$. By  Lemma \ref{rig-lemma: lemma} we have 
$$\hat{\mu}^{\sigma,\xi}_{\tilde{v}^l}(U) \le \liminf_{k \to \infty} \hat{\mu}^{\sigma,\xi}_{\tilde{v}^l_k}(U) $$
for all $\sigma\ge 0$, $\xi \in S^1$ and for every open set $U \subset  \Omega$. Let $\kappa_1 = \int_{S^1} |\nu \cdot \xi| \,d{\cal H}^1(\xi)$ for some $\nu \in S^1$ which clearly does not depend on the particular choice of $\nu$. Using Fatou's lemma and \eqref{rig-eq: lemma2} we compute  
\begin{align*}
\liminf_{k\to \infty} {\cal H}^1(J_{y_k}) &+ C\rho_l \ge \liminf_{k \to \infty}  \int_{J^l_k} \theta_\sigma(|[\tilde{v}^l_k]|) \, d{\cal H}^1 \\
& \ge \kappa_1^{-1} \int_{S^1} \liminf_{k \to \infty} \int_{J^l_k} \theta_\sigma(|[\tilde{v}^l_k](x)|) |\nu_{\tilde{v}^l_k}(x) \cdot \xi| \, d{\cal H}^1(x) \, d{\cal H}^1(\xi) \\
& \ge \kappa_1^{-1} \int_{S^1} \liminf_{k \to \infty} \hat{\mu}^{\sigma,\xi}_{\tilde{v}^l_k}(\Omega_{\rho_l}) \, d{\cal H}^1(\xi) \ge \kappa_1^{-1} \int_{S^1} \hat{\mu}^{\sigma,\xi}_{\tilde{v}^l}(\Omega_{\rho_l}) \, d{\cal H}^1(\xi). 
\end{align*}
We pass to the limit $l \to \infty$ (i.e. $\rho_l \to 0$) and obtain \BBB by Lemma \ref{rig-lemma: lemma} and \EEE the dominated convergence theorem
$$\liminf_{k\to \infty} {\cal H}^1(J_{y_k}) \ge \kappa_1^{-1} \int_{S^1} \hat{\mu}^{\sigma,\xi}_{\tilde{v}}(\Omega) \, d{\cal H}^1(\xi).$$
Recall that $\theta_\sigma \to 1$ pointwise for $\sigma \to 0$. Now letting $\sigma \to 0$ we obtain by the dominated convergence theorem and \eqref{rig-eq: lemma2}
\begin{align*}
\liminf_{k \to \infty} {\cal H}^1(J_{y_k}) & \ge \kappa_1^{-1} \int_{S^1} \hat{\mu}^{0,\xi}_{\tilde{v}}(\Omega) \, d{\cal H}^1(\xi) \\
& = \kappa_1^{-1} \int_{S^1} \int_{J^\xi_{\tilde{v}}} |\nu_{\tilde{v}}(x) \cdot \xi| \, d{\cal H}^1(x)\, d{\cal H}^1(\xi) =  {\cal H}^1( J_{\tilde{v}} ).
\end{align*}
This gives \eqref{rig-eq: compa1} and completes the proof. \eop

\begin{rem}\label{rem: NNNN}
{\normalfont

Using \eqref{rig-eq: comp2}(ii), \eqref{rig-eq: comp3}(i) and arguing as in \eqref{rig-eq: ff}, \eqref{rig-eq: ff2}, we observe that all sequences  $(u_k, \mathcal{P}_k, T_k)$  in Definition \ref{def:conv} satisfy $\Vert e(\nabla T_k^T \nabla u_k) \Vert_{L^2(\Omega)} \le C$ for $C>0$ only depending on $\sup_k E_{\eps_k}(y_k)$, $\Omega$, and the constant in \eqref{eq:W}.} 
\end{rem}

At the end of this section we briefly note that our compactness result provides   an alternative proof of the piecewise rigidity result given in Theorem \ref{rig-cor: cgp} (at least in a planar setting).

\noindent {\em Proof of Theorem  \ref{rig-cor: cgp} for $d=2$.}  Let $y \in SBV(\Omega)$ with ${\cal H}^1(J_y) < \infty$ as well as $\int_\Omega \dist^2(\nabla y,SO(2)) = 0$ be given. \BBB First, assume $y \in L^\infty(\Omega)$. \EEE Define an arbitrary infinitesimal sequence $(\eps_k)_k$ and the sequence $y_k = y$ for all $k \in \N$. Applying Theorem \ref{rig-th: comp1} we obtain piecewise rigid motions $T,T_k$ such that $T_k \to T$,  $\nabla T_k \to \nabla T$ in $L^2(\Omega)$ by \eqref{rig-eq: comp1} up to passing to a subsequence. Moreover, $y_k - T_k \to 0$ a.e.  in $\Omega$ for $k \to \infty$ by \eqref{rig-eq: comp2}(i). This implies $y = T$ is a piecewise rigid motion. \BBB If $y \notin L^\infty(\Omega)$, using the $BV$ coarea formula we can approximate $y$  by a  sequence $y \chi_{\Phi_k} + \id\chi_{\Omega \setminus \Phi_k} \in SBV(\Omega) \cap L^\infty(\Omega)$ with $\sup_{k} {\cal H}^1(\partial^* \Phi_k )<\infty$, $|\Phi_k| \to 0$ for $k \to \infty$ and conclude by Theorem \ref{th: piecewise const}. \EEE \eop

%----------------------------------------------------------------------------------------------------------------------
\section{Admissible  \BBB and \EEE coarsest partitions and limiting configurations}\label{rig-sec: sub, comp2}
%-----------------------------------------------------------------------------------------------------------------------

In this section we will prove Theorem \ref{rig-th: comp2}. Let $(y_k)_k$ be a (sub-)sequence as considered in Theorem \ref{rig-th: comp1}. Recall Definition \ref{rig-def: ad,coar}. For notational convenience we will drop the  dependence of $(y_k)_k$ in the sets ${\cal Z}_P, {\cal Z}_u, {\cal Z}_T$. We introduce a partial order on the admissible partitions ${\cal Z}_P$: Given two partitions ${\cal P}^1 :=(P^1_j)_j, {\cal P}^2:=(P^2_j)_j$ in ${\cal Z}_P$ we say ${\cal P}^2 \ge {\cal P}^1$ if  \BBB 
\begin{align}\label{new10}
\text{for all \  $P^1_{j_1}$ \  there exists \ $P^2_{j_2}$ \  such that \ $|P^1_{j_1} \setminus P^2_{j_2}|=0$.}
\end{align}
Note that Theorem \ref{th: local structure} implies $\bigcup_j \partial^* P^1_j \supset \bigcup_j \partial^* P^2_j$  \EEE up to an ${\cal H}^1$-negligible set.  We observe that if  ${\cal P}^1 \ge {\cal P}^2$ and ${\cal P}^2 \ge {\cal P}^1$, abbreviated by ${\cal P}^1 = {\cal P}^2$ hereafter, then the Caccioppoli partitions coincide: After a possible reordering of the sets we find $|P^1_j \triangle P^2_j| = 0$ for all $j \in \N$.

We begin with the observation that the piecewise rigid motion is uniquely determined in the limit.

\begin{lemma}\label{rig-lemma: comp2.5}
Let $(y_k)_k$ be a (sub-)sequence as considered in Theorem \ref{rig-th: comp1}. Then there is a unique $T \in {\cal Z}_T$.
\end{lemma}

\Proof Assume there are $T,\hat{T}  \in {\cal Z}_T$. Let  $(u,{\cal P},T), (\hat{u},\hat{\cal P},\hat{T}) \in {\cal D}_\infty$
according to Definition \ref{rig-def: ad,coar}(ii) and \BBB let $(u_k,{\cal P}_k,T_k), (\hat{u}_k,\hat{\cal P}_k,\hat{T}_k) \in {\cal D}$ for $k \in \N$ be triples given by Definition \ref{def:conv}. \EEE As $u_k - \hat{u}_k - (\eps_k^{-1/2} (T_k - \hat{T}_k)) \to 0$ a.e. by \eqref{rig-eq: comp2}(i) and $u_k - \hat{u}_k$ converges pointwise a.e. (and the limits lie in $\R$ a.e.) by \eqref{rig-eq: comp1-2}(i), we get $T_k - \hat{T}_k   \to 0$ pointwise almost everywhere. \BBB By \eqref{rig-eq: comp1}(ii) \EEE this implies $T = \hat{T}$. \eop

From now on $T$ will always stand for the rigid motion given by Lemma \ref{rig-lemma: comp2.5}.

\subsection{Equivalent characterization of the coarsest partition} 

We state a lemma giving an  equivalent characterization of the coarsest partition  (recall Definition \ref{rig-def: ad,coar}(iv)).

\begin{lemma}\label{rig-lemma: comp2.1}
Let $(y_k)_k$ be a (sub-)sequence as considered in Theorem \ref{rig-th: comp1}. Then ${\cal P} \in {\cal Z}_P$ is coarsest if and only if it is a maximal element in the partial order $({\cal Z}_P, \ge)$, i.e. $\hat{\cal P} \ge {\cal P} $ implies $\hat{\cal P} = {\cal P} $.
\end{lemma}

\Proof (1) Assume ${\cal P}  = (P_j)_j$ was not coarsest. According to Definition \ref{rig-def: ad,coar}(iv)  let $u$ and $(u_k, {\cal P}_k, T_k)\in {\cal D}$   be given such that $(u,{\cal P},T) \in {\cal D}_\infty$ and \eqref{rig-eq: comp2}-\eqref{rig-eq: comp3} hold. Without restriction, possibly passing to a subsequence \BBB and reordering the partition, \EEE we assume that $  \eps_k^{-1/2} \big(|R^k_1 - R^k_2| + |b_1^k - b_2^k| \big) \le C$ for all $k \in \N$  (cf. \eqref{rig-eq: toinfty}). By \eqref{rig-eq: linearization} we obtain $A^k \in \R^{2 \times 2}_{\rm skew}$ for $k \in \N$ with $|A^k| \le C$ such that $R_2^k - R^k_1 = R_1^k ((R_1^k)^T R^k_2 - \Id)  = R_1^k (\sqrt{\eps_k} A^k + O(\eps_k))$. Passing to a  (not relabeled) subsequence we then obtain for all $x \in \Omega$
\begin{align}\label{rig-eq: Tdef}
\begin{split}
S(x) & := \lim_{k \to \infty} \frac{1}{\sqrt{\eps_k}} \big((R^k_2 - R^k_1)\,x + b_2^k - b_1^k \big) \\
& = \lim_{k \to \infty} \frac{1}{\sqrt{\eps_k}} \big( \sqrt{\eps_k} R^k_1 A^k\,x + b_2^k - b_1^k \big) + O(\sqrt{\eps_k}) = R A\,x +b
\end{split}
\end{align}
for some $A \in \R^{2 \times 2}_{\rm skew}$, $b \in \R^2$ and $R = \lim_{k \to \infty} R^k_1$. We now introduce $\hat{ \cal P}_k$, $\hat{ \cal P}$, $\hat{T}_k$, $\hat{u}_k, \hat{u}$ as follows. Let $\hat{P}^k_1 = P^k_1 \cup P^k_2$, $\hat{P}^k_2 = \emptyset$, $\hat{P}^k_j = P^k_j$ for $j \ge 3$ and likewise for the limiting partition $\hat{\cal P}$. Let $\hat{T}_k(x) = R^k_1 \, x + b^k_1$ for $x \in \hat{P}^k_1$ and $\hat{T}_k(x) = T_k(x)$ for $x \in \Omega \setminus \hat{P}^k_1$. Furthermore, we let 
$$\hat{u}_k = u_k + \frac{1}{\sqrt{\eps_k}} \big((R^k_2 - R^k_1)\cdot + b_2^k - b_1^k \big) \chi_{P^k_2}$$
and $\hat{u} = u + (RA\cdot+b) \chi_{P_2}$   (see \eqref{rig-eq: Tdef}).  \BBB We now show that $(\hat{u}_k, \hat{\cal P}_k, \hat{T}_k)$ converges to $(\hat{u},\hat{\cal P},T)$ in the sense of \eqref{rig-eq: comp2}-\eqref{rig-eq: comp3}. First, \EEE \eqref{rig-eq: comp2}(i) clearly holds since $\hat{T}_k - T_k = \big((R^k_1 - R^k_2)\cdot+ b_1^k - b_2^k \big)\chi_{P_2^k}$. Moreover, we derive that \eqref{rig-eq: comp1} holds as $|R^k_1 - R^k_2| + |b_1^k - b_2^k| \to 0$ for $k \to \infty$. Since  ${\cal H}^1(J_u \setminus \bigcup_j \partial^* P_j) + {\cal H}^1(\bigcup_j \partial^* P_j \cap \Omega)  \ge {\cal H}^1(J_{\hat{u}} \setminus \bigcup_j \partial^* \hat{P}_j) + {\cal H}^1(\bigcup_j \partial^* \hat{P}_j \cap \Omega)$, also \eqref{rig-eq: comp3}(ii) is satisfied.

\BBB As $|R^k_2 - R^k_1| \le C\sqrt{\eps_k}$,   we note that $\Vert \nabla \hat{u}_k\Vert_{L^\infty(\Omega)} \le c\eps_k^{-1/8} $ for $c>0$ large enough and thus \eqref{rig-eq: comp2}(ii) holds. \EEE  It remains to verify \eqref{rig-eq: comp1-2}  and \eqref{rig-eq: comp3}(i). First, \eqref{rig-eq: comp1-2}(i) follows from \eqref{rig-eq: Tdef} and the definition of $\hat{u}$.  We use $R^k_2 = R_1^k + \sqrt{\eps_k} R^k_1 A^k +  O(\eps_k)$, $|A^k|\le C$ and $\Vert \nabla u_k \Vert_{L^\infty(\Omega)}\le c\eps_k^{-1/8}$  to find a.e. on $P_2^k$
\begin{align}\label{rig-eq: bar-nonbar}
\nabla \hat{T}_k^T \nabla \hat{u}_k =  (R_1^k)^T \nabla u_k + A^k + O(\sqrt{\eps_k}) = (R_2^k)^T \nabla u_k + A^k + O(\eps^{3/8}_k).
\end{align}
 Now we get 
\begin{align*}
\chi_{\hat{P}^k_1}  e(\nabla \hat{T}_k^T \nabla \hat{u}_k) & = \sum\nolimits_{j=1,2}\chi_{P^k_j}  e((R^k_j)^T \nabla u_k) + \chi_{P^k_2}  e( A^k) + \BBB O(\eps^{3/8}_k) \EEE \\
& \rightharpoonup \sum\nolimits_{j=1,2} \chi_{P_j} e(R_j^T \nabla u) = \chi_{\hat{P}_1} e(R^T\nabla \hat{u}) \BBB = \chi_{\hat{P}_1} e(\nabla T^T\nabla \hat{u}) \EEE
\end{align*}
weakly in $L^2(\Omega,\R^{2 \times 2}_{\rm sym}).$ \BBB  This gives \eqref{rig-eq: comp1-2}(ii). By the assumptions on $W$, a Taylor expansion yields \EEE   $W(G) = \frac{1}{2} Q( e(G-\Id)) + \omega_{\rm \BBB W \EEE}(G-\Id)$ for $G \in  \R^{2 \times 2}$, where $\sup\lbrace |F|^{-3}\omega_W(F): |F| \le 1\rbrace \le C$ and $Q  =  D^2 W(\Id)$. Thus, we obtain by \eqref{rig-eq: bar-nonbar} and \BBB $\Vert \nabla u_k \Vert_{L^\infty(\Omega)}\le c\eps_k^{-1/8}$ \EEE
\begin{align*}
\frac{1}{\eps_k} \int_{P_2^k}  \hspace{-0.1cm} W(\Id + \sqrt{\eps_k} \nabla \hat{T}_k^T \nabla \hat{u}_k)  &  = \int_{P_2^k}  \Big(\frac{1}{2}Q(e(\nabla \hat{T}_k^T \nabla \hat{u}_k)) + \frac{1}{\eps_k}\omega_W(\sqrt{\eps_k}\nabla\hat{T}_k^T \nabla \hat{u}_k) \Big) \\
& = \int_{P_2^k} \hspace{-0.1cm} \Big(\frac{1}{2}Q(e(\nabla T_k^T \nabla u_k)) +  \frac{\omega_W(\sqrt{\eps_k} \nabla \hat{u}_k)}{\eps_k}\Big) + O(\eps_k^{\frac{\BBB 1 \EEE}{4}})
\end{align*}
and likewise
\begin{align*}
\frac{1}{\eps_k} \int_{P_2^k}  W (\Id + \sqrt{\eps_k} \nabla T_k^T \nabla u_k) = \int_{P_2^k}    \Big(\frac{1}{2}Q(e(\nabla T_k^T \nabla u_k)) +  \frac{1}{\eps_k}\omega_W(\sqrt{\eps_k}  \nabla u_k)\Big). 
\end{align*}
In both estimates the second terms converge to $0$ using $\Vert \nabla u_k \Vert_\infty + \Vert \nabla \hat{u}_k \Vert_\infty \le c\eps_k^{-1/8} $ and arguing as in \eqref{rig-eq: ff2}. Consequently,  we get  
\begin{align}\label{rig-eq: bar-nonbar2} 
\frac{1}{\eps_k} \int_{P_2^k}   W(\Id + \sqrt{\eps_k} \nabla \hat{T}_k^T \nabla \hat{u}_k)  = \frac{1}{\eps_k} \int_{P_2^k}   W(\Id + \sqrt{\eps_k} \nabla T_k^T \nabla u_k) + o(1)
\end{align}
for $\eps_k \to 0$, i.e.  \eqref{rig-eq: comp3}(i) holds.  Therefore,   $\hat{\cal P}$ is an admissible partition and thus ${\cal P}$ is not maximal. 

\smallskip

(2) Conversely, assume that ${\cal P} = (P_j)_j$ was not maximal, i.e. we find $\hat{\cal P} = (\hat{P}_j)_j$ with $\hat{\cal P} \ge {\cal P}$, $\hat{\cal P} \neq {\cal P}$. \BBB Upon reordering \EEE we may assume  that $ P_1 \cap \hat{P}_1$ and $P_2 \cap \hat{P}_1$ have positive ${\cal L}^2$-measure.  According to Definition \ref{rig-def: ad,coar}(i)  let $u,\hat{u}$ and $(u_k, {\cal P}_k,T_k), (\hat{u}_k, \hat{\cal P}_k, \hat{T}_k) \in {\cal D}$  be given such that $(u,{\cal P},T), (\hat{u},\hat{\cal P},T) \in {\cal D}_\infty$ and \eqref{rig-eq: comp2}-\eqref{rig-eq: comp3} hold. As by \eqref{rig-eq: comp1-2}(i) $u_k$ and $\hat{u}_k$  convergence pointwise a.e., by \eqref{rig-eq: comp2} also $  \eps_k^{-1/2}(T_k - \hat{T}_k)$ converges pointwise a.e. (and the limits  lie in $\R$ a.e.). But this implies $  \eps_k^{-1/2} \big(|R^k_j - \hat{R}^k_1| + |b_j^k - \hat{b}_1^k| \big) \le C$ for $j=1,2$ and $k \in \N$.  \BBB Then  $  \eps_k^{-1/2}\big(|R^k_1 - {R}^k_2| + |b_1^k - {b}_2^k| \big) \le C$ by the triangle inequality. Thus, \EEE \eqref{rig-eq: toinfty} is violated and ${\cal P}$ is not a coarsest partition. \eop

 The alternative characterization now directly implies that there is at most one coarsest partition.

\begin{lemma}\label{rig-lemma: comp2.2}
Let $(y_k)_k$ be a (sub-)sequence as considered in Theorem \ref{rig-th: comp1}. Then there is at most one maximal element in $({\cal Z}_P, \ge)$.
\end{lemma}

\Proof Assume there are two maximal elements ${\cal P}^1 =(P^1_j)_j,{\cal P}^2 = (P^2_j)_j \in {\cal Z}_P$ with ${\cal P}^1 \neq {\cal P}^2$. As before, without restriction we may assume that $P^1_1 \cap P^2_1$ and $P^1_2 \cap P^2_1$ have positive ${\cal L}^2$-measure. We proceed as in the proof of Lemma \ref{rig-lemma: comp2.1}(2) to see that   ${\cal P}^1$ is not coarsest and thus not a maximal element in $({\cal Z}_P, \ge)$.  \eop

\subsection{Admissible configurations}

We now analyze the admissible configurations if the partitions are given.  Recall the definition of the set  of piecewise infinitesimal rigid motions ${\cal A}({\cal P})$ below \eqref{rig-eq: defA}.

\begin{lemma}\label{rig-lemma: comp2.4}
Let $(y_k)_k$ be a (sub-)sequence as considered in Theorem \ref{rig-th: comp1} and \BBB let  $T \in {\cal Z}_T$ be the unique mapping given by Lemma \ref{rig-lemma: comp2.5}. \EEE Let ${\cal P},\hat{\cal P} \in {\cal Z}_P$ such that $\hat{\cal P} \ge {\cal P}$ and  $\hat{u} \in {\cal Z}_u(\hat{\cal P})$. Then ${\cal Z}_u({\cal P}) = \hat{u} + \nabla T {\cal A}({\cal P})$.
\end{lemma}

\Proof (1) To see  ${\cal Z}_u({\cal P}) \subset \hat{u} +  \nabla T {\cal A}({\cal P})$, we have to show that $u - \hat{u} \in  \nabla T {\cal A}({\cal P})$ for all $u \in {\cal Z}_u({\cal P})$.  To this end, consider $P_{j} \in {\cal P}$, $\hat{P}_{i} \in \hat{\cal P}$ such that  $|P_{j} \setminus \hat{P}_{i}| = 0$. Let $(u_k,{\cal P}_k,T_k), (\hat{u}_k,\hat{\cal P}_k,\hat{T}_k) \in {\cal D}$   be given according to \BBB Definition \ref{def:conv}. \EEE As $u_k - \hat{u}_k$ and thus $  \eps_k^{-1/2}(T_k - \hat{T}_k)$ converge pointwise a.e. \BBB by \eqref{rig-eq: comp2}(i), \eqref{rig-eq: comp1-2}(i),  \EEE  we find $|R^k_{j} - \hat{R}^k_{i}| + |b^k_{j} - \hat{b}^k_{i}| \le C\sqrt{\eps_k}$. Repeating the argument in \eqref{rig-eq: Tdef} we find some $A \in \R^{2 \times 2}_{\rm skew}$, $b \in \R^2$ such that for a.e. $x \in P_j$ 
$$u(x)- \hat{u}(x) = \lim_{k \to \infty} u_k(x)- \hat{u}_k(x) = \lim_{k \to \infty} \eps_k^{-1/2} (\hat{T}_k(x)- {T}_k(x)) =  \nabla T(x) (A \,x + b).$$

(2) Conversely, to see  ${\cal Z}_u({\cal P}) \supset \hat{u} + \nabla T{\cal A}({\cal P})$ we first consider the special case ${\cal P} = \hat{\cal P} = (P_j)_j$. Let $\bar{u} \in {\cal Z}_u( { \cal P})$ and $\bar{A} = \sum_{j} (A_j\cdot + d_j) \chi_{P_j} \BBB \in {\cal A}({\cal P}) \EEE $ be given. We have to show that $u := \bar{u}+ \nabla T \bar{A} \in {\cal Z}_u( {\cal P})$.

According to Definition \ref{rig-def: ad,coar}(iii) let $(\bar{u}_k, {\cal P}_k, \bar{T}_k) \in {\cal D}$ be given such that \eqref{rig-eq: comp2}-\eqref{rig-eq: comp3} hold \BBB with the limiting triple $(\bar{u},{\cal P},T)$. \EEE   Assume that $\bar{T}_k$ has the form $\bar{T}_k = \sum_j( \bar{R}^k_j \cdot + \bar{b}^k_j) \chi_{P^k_j}$. Now choose  $R^k_j$ such that  $| R^k_j - \bar{R}^k_j (\Id - \sqrt{\eps_k} A_{j})| = \dist(\bar{R}^k_j (\Id - \sqrt{\eps_k} A_{j}), SO(2))$ and let $b^k_j = \bar{b}^k_j -\sqrt{\eps_k} \bar{R}^k_j  d_j$. \BBB By \eqref{rig-eq: linearization} we have 
\begin{align}\label{new5}
R^k_j = \bar{R}^k_j - \sqrt{\eps_k}  \bar{R}^k_j  A_{j}- \omega_{j,k} \text{ with   $|\omega_{j,k}| \le C\eps_k|A_j|^2$   for all $j \in \N$.}
\end{align}
Let  $I_k = \lbrace j \in \N: |A_j| + |d_j| \le \eps_k^{-1/8} \rbrace$ and $V_k = \bigcup_{j \in \N \setminus I_k} P_j^k$. Note that $|V_k| \to 0$  for $k \to \infty$  and $|\eps_k^{-1/2}\omega_{j,k}| \le C\eps^{1/4}_k$ for $j \in I_k$. \EEE  Define
$$ \BBB T_k = \sum\nolimits_{j \in I_k}  (R^k_j \cdot + b^k_j ) \chi_{{P}^k_j} +  \bar{T}_k\chi_{V_k}, \EEE \ \ \ \  u_k = \bar{u}_k + \frac{1}{\sqrt{\eps_k}}(\bar{T}_k - T_k). $$\BBB We now show that $(u_k, {\cal P}_k,T_k)$ converges to $(u,{\cal P},T)$ in the sense of \eqref{rig-eq: comp2}-\eqref{rig-eq: comp3}. First, \EEE \eqref{rig-eq: comp2}(i) clearly holds. Moreover, ${\cal H}^1(J_u \setminus \bigcup_j \partial^* P_j) = {\cal H}^1(J_{\bar{u}} \setminus \bigcup_j \partial^* P_j)$ and thus \eqref{rig-eq: comp3}(ii) is satisfied. \BBB By \eqref{new5} \EEE    we find 
$$ T_k =  \bar{T}_k -  \sum\nolimits_{\BBB j\in I_k \EEE } \big(\sqrt{\eps_k} \bar{R}^k_j A_j \cdot + \omega_{j,k} \cdot + \sqrt{\eps_k}\bar{R}_j^k d_{j}\big) \chi_{{P}^k_j} \to T$$
in measure for $k \to \infty$. Then it is not hard to see that  $T_k \to T$ and $\nabla T_k \to \nabla T$ in   $L^2(\Omega)$  which gives  \eqref{rig-eq: comp1}. Likewise, we obtain
\begin{align*}
 u_k - \bar{u}_k   &= \frac{1}{\sqrt{\eps_k}}\big(\bar{T}_k - T_k\big) = \sum\nolimits_{j \in I_k} \Big( \bar{R}^k_j  (A_j\cdot +  d_j)  + \frac{1}{\sqrt{\eps_k}}\omega_{j,k}\cdot\Big)\chi_{{P}^k_j} \\ & \to \nabla T \sum\nolimits_j (A_j \cdot + d_j)\chi_{P_j}  = \nabla T \bar{A}
\end{align*}
pointwise a.e. which   implies  $u_k \to \bar{u} + \nabla T \bar{A}$ and shows \eqref{rig-eq: comp1-2}(i). \BBB   By the definition of $V_k$ we have \EEE   
$$\Vert \nabla u_k - \nabla \bar{u}_k\Vert_{L^\infty(\BBB \Omega \EEE)} \le \Vert \sum\nolimits_j \chi_{P^k_j}\big( \BBB \bar{R}_j^k A_j \EEE + \eps_k^{-1/2}\omega_{j,k}\big)\Vert_{L^\infty(\Omega \setminus V_k)}\le C\eps_k^{-1/8}.$$Therefore, \BBB $\Vert \nabla u_k \Vert_\infty \le c\eps_k^{-1/8}$ for $c$ large enough, which shows \eqref{rig-eq: comp2}(ii). \EEE    Arguing  as in \eqref{rig-eq: bar-nonbar} we get by $\Vert \nabla \bar{u}_k \Vert_{L^\infty(\Omega)} \le c\eps_k^{-1/8}$ and \BBB  \eqref{new5}\EEE
\begin{align*}
(R_j^k)^T \nabla u_k(x) &=   (R_j^k)^T \nabla \bar{u}_k(x) + \BBB (R_j^k)^T\bar{R}_j^k \EEE A_j + (R_j^k)^T \eps_k^{-1/2}w_{j,k}\\
& =(\bar{R}_j^k)^T \nabla \bar{u}_k(x) + A_j + O(\eps_k^{1/4})  
\end{align*}
for  a.e. $x \in P_j^k$, \BBB $j \in I_k$. \EEE  Thus,  \eqref{rig-eq: comp1-2}(ii) follows from  the fact that \eqref{rig-eq: comp1-2}(ii) holds for the sequence $\bar{u}_k$ and  
\begin{align*}
\sum\nolimits_{\BBB j\in I_k}  \int\nolimits_{{P}^k_j} |e\big(( R^k_j)^T\nabla u_k \big)& -  e\big(( \bar{R}^k_j)^T\nabla \bar{u}_k \big)|^2 \le  C\eps_k^{1/2} \to 0.
\end{align*}
Finally, the above estimates together with a similar argumentation as in   \eqref{rig-eq: bar-nonbar2}  yield \eqref{rig-eq: comp3}(i). 

In the general case we have to show $u := \hat{u}+ \nabla T \bar{A} \in {\cal Z}_u( {\cal P})$ for given $\hat{u} \in {\cal Z}_u( {\cal \hat{P}})$, $\hat{\cal P} \ge {\cal P}$, and $\bar{A} \in {\cal A}({\cal P})$. As ${\cal P} \in {\cal Z}_P$, we find some $\bar{u} \in {\cal Z}_u({\cal P})$ which by (1) satisfies $\bar{u} - \hat{u} = \nabla T \hat{A}$ for some $\hat{A} \in {\cal A}({\cal P})$. Thus, we get $u = \bar{u} + \nabla T(\bar{A} - \hat{A})$ and by the special case in (2) we know that $u \in {\cal Z}_u({\cal P})$, as desired. \eop

%We observe that for each $i \ge i_0$ we have $\bigcup_j \partial^* P_{i,j} \subset \bigcup_j \partial^* P_{i_0,j}$ (up to an ${\cal H}^1$-negligible set). Thus, for each $k \in \N$ there are (coarsened) partitions ${\cal P}^k_i = (P^k_{i,j})_j$ with $\bigcup_j \partial^* P^k_{i,j} = \bigcup_j \partial^* P_{i,j} \setminus \big(\bigcup_{j\ge k} \partial^* P_{i_0,j} \setminus \partial^* (\bigcup_{j \ge k} P_{i_0,j})\big)$ up to ${\cal H}^1$-negliglible sets for all $i \ge i_0$. Observe that typically ${\cal P}^k_i$ are not elements of $\lbrace {\cal P}_i: i \in I \rbrace$, but satisfy 
% $${\cal H}^1\Big(\bigcup\nolimits_j\partial^* P_{i,j} \setminus \bigcup\nolimits_j\partial^* P^k_{i,j}\Big) \le \omega(k)$$
%After identifying partitions whose boundaries only differ by ${\cal H}^1$-negligible sets we find that each 

\subsection{Existence of coarsest partitions}

To guarantee existence of coarsest partitions we show that each totally ordered subset has upper bounds such that afterwards we may apply Zorn's lemma.

\begin{lemma}\label{rig-lemma: comp2.3}
Let $(y_k)_k$ be a (sub-)sequence as considered in Theorem \ref{rig-th: comp1}.  Let $I$ be an arbitrary index set and let $\lbrace {\cal P}_i = (P_{i,j})_j: i \in I\rbrace \subset{\cal Z}_P$ be a totally ordered subset, i.e. for each $i_1,i_2 \in I$ we have ${\cal P}_{i_1} \le {\cal P}_{i_2}$ or ${\cal P}_{i_2} \le {\cal P}_{i_1}$. Then there is a partition ${\cal P} \in {\cal Z}_P$ with ${\cal P}_i \le {\cal P}$ for all $i \in I$. 
\end{lemma}

\Proof \BBB \emph{Step I:} \EEE To prove the existence of an upper bound we first show that it suffices to consider a suitable countable subset of  $\lbrace{\cal P}_i: i \in I \rbrace$. For notational convenience we write $i_1 \le i_2$ for $i_1, i_2 \in I$ if ${\cal P}_{i_1} \le {\cal P}_{i_2}$. Choose an arbitrary $i_0 \in I$ and note that it suffices to find an upper bound for all $i \ge i_0$. \BBB   For  each $k \in \N$ we introduce partitions ${\cal P}^k_i = (P^k_{i,j})_{j \ge 0}$ consisting of the components $P^k_{i,j} = P_{i,j} \setminus \bigcup\nolimits_{l\ge k}P_{i_0,l}$ for $j \in \N$ and $P^k_{i,0} = \bigcup\nolimits_{l\ge k}P_{i_0,l}$.  (Note that the partitions ${\cal P}^k_i$ are possibly not ordered.) By \eqref{new10}  we get that ${\cal P}^k_{i_1} \le {\cal P}^k_{i_2}$ if $i_0 \le i_1\le i_2$. Typically, ${\cal P}^k_i$ are not elements of $\lbrace {\cal P}_i: i \in I \rbrace$, but satisfy for $i \ge i_0$ 
\begin{align}\label{new11}
\begin{split}
&\text{$|P_{i,j} \triangle P^k_{i,j}| \le  \big|\bigcup\nolimits_{l\ge k}P_{i_0,l} \big| \le \omega(k)$ for all $j \ge 0$}
\end{split}
\end{align}
with $\omega(k) \to 0$ for $k \to \infty$,  where we set $P_{i,0} = \emptyset$. For all $k \in \N$ we observe that $\lbrace {\cal P}^k_i: i \ge i_0\rbrace$ contains only a finite number of different elements and therefore contains a maximal element ${\cal P}^k = (P^k_j)_j$. Now we can choose $i_0 \le i_1 \le i_2 \le \ldots$ such that ${\cal P}^k = {\cal P}^k_{i_k}$ for $k\in\N$.  It now suffices to construct an upper bound ${\cal P} = (P_j)_j \in {\cal Z}_P$ with ${\cal P} \ge {\cal P}_{i_k}$ for all $k \in \N$. Indeed, we then obtain ${\cal P} \ge {\cal P}_{i}$ for all $i_0 \le i$ as follows:

For each $P_{i,j}$ and each $k \in \N$ we find $P_{j_k}$ with $|P_{i,j} \setminus P_{j_k}| \le 2\omega(k)$. In fact, using repetitively \eqref{new10} and \eqref{new11} we get $j',j_k$ such that $|P_{i,j} \setminus P_{i,j}^k| \le \omega(k)$, $|P_{i,j}^k \setminus P^k_{i_k,j'}|=0$, $|P^k_{i_k,j'} \setminus P_{i_k,j'}| \le \omega(k)$, $|P_{i_k,j'} \setminus P_{j_k}|=0$ and thus $|P_{i,j} \setminus P_{j_k}| \le 2\omega(k)$. As $\omega(k) \to 0$ for $k \to \infty$ and ${\cal P}$ contains only a finite number of components with $\mathcal{L}^2$-measure larger than $\frac{1}{2}|P_{i,j}|$, we indeed find $P_{j_*}$ with $|P_{i,j} \setminus P_{j_*}| = 0$, as desired. \EEE

\smallskip

\BBB \emph{Step II:} \EEE Now consider the totally ordered sequence of partitions $({\cal P}_{i_k})_k$. For notational convenience we will denote the sequence by $({\cal P}_{i})_{i \in \N}$ in the following. By the compactness theorem for Caccioppoli partitions (see Theorem \ref{th: comp cacciop}) we get an (ordered) Caccioppoli partition ${\cal P}=(P_j)_j$ such that \BBB $|P_{i,j} \triangle P_j| \to 0$ \EEE for $i \to \infty$ for all $j \in \N$.  Thus, for all $j \in \N$ there exists $I_j \in \N$ such that $|P_{i_1,j} \setminus P_{i_2,j}| \le \frac{1}{2}|P_{i_1,j}|$ for all $I_j \le i_1 \le i_2$.  As  $({\cal P}_{i})_{i \in \N}$ is totally ordered, \eqref{new10} then gives $|P_{i_1,j} \setminus P_{i_2,j}|=0$ for all $I_j \le i_1 \le i_2$ and this monotonicity yields $|P_{i_1,j} \setminus P_j|=0$ for $i_1 \ge I_j$. Eventually, fixing $P_{i,j}$ for $i,j \in \N$, by the above arguments there exists $j' \in \N$ such that  $|P_{i,j} \setminus P_{i',j'}|=0$ for all $i'$ large enough and thus $|P_{i,j} \setminus P_{j'}|=0$.  

This implies ${\cal P} \ge {\cal P}_i$ for all $i \in \N$ and therefore it suffices to show that ${\cal P} \in {\cal Z}_P$. To this end, we will construct partitions ${\cal P}^n$, rigid motions $T_n \in {\cal R}({\cal P}^n)$ and a limiting function $u$ by a diagonal sequence argument. 

For all $i \in \N$, according to Definition \ref{rig-def: ad,coar}(i), we find $(u^k_i, {\cal P}^k_i, T^k_i) \in {\cal D}$, an admissible limiting configurations $u_i \in {\cal Z}_u( {\cal P}_i)$ and \BBB $T \in {\cal Z}_T$ as in Lemma \ref{rig-lemma: comp2.5} \EEE such that \eqref{rig-eq: comp2}-\eqref{rig-eq: comp3} hold  as $k \to \infty$. The strategy is to select $u_i$ in a suitable way such that we find a limiting configuration $u \in GSBD(\Omega)$ with
\begin{align}\label{rig-eq: vi to v}
\begin{split}
&u_i \to u \   \text{ a.e.  in $\Omega$, } \\ &  e(\nabla T^T \nabla u_i) \rightharpoonup e(\nabla T^T \nabla u) \ \ \ \text{weakly in} \ L^2(\Omega,\R^{2\times 2}_{\rm sym}),\\ & \liminf\nolimits_{i \to \infty} {\cal H}^1(J_{u_i}) \ge {\cal H}^1(J_u).
\end{split}
\end{align}
\BBB We defer the selection of the sequence $(u_i)_i$ to Step III below. \EEE Then we can choose a diagonal sequence $(\bar{u}_n) := (u^{k(n)}_n)_n$ converging to the triple $(u,{\cal P},T)$ in the sense of \eqref{rig-eq: comp2}-\eqref{rig-eq: comp3}. Indeed, $k(n)$ can be selected such that letting $\bar{\cal P}^n =  (\bar{P}^n_j)_j =  {\cal P}^{k(n)}_n$ and $\bar{T}_n = T^{k(n)}_n \in {\cal R}(\bar{\cal P}^n)$,  we find \BBB $|\bar{P}^n_j \triangle P_j| \to 0$ \EEE for all $j \in \N$ (even $\sum_j|\bar{P}^n_j \triangle P_j| \to 0$, cf. below Theorem \ref{th: comp cacciop}) and  $\bar{T}_n \to T$, $\nabla \bar{T}_n \to \nabla T$ in $L^2(\Omega)$. This gives \eqref{rig-eq: comp1}.  Moreover, as measure convergence is metrizable, this can be done in a way  that $\bar{u}_n\to u$ in measure and $\bar{u}_n - \eps_n^{-1/2}(y_n - \bar{T}_n) \to 0$ in measure. Then, possibly passing to a further subsequence, we can assume that the convergence also holds a.e. in $\Omega$ and thus \eqref{rig-eq: comp2}(i), \eqref{rig-eq: comp1-2}(i) are satisfied. 

 Likewise, \eqref{rig-eq: comp1-2}(ii) can be achieved by \eqref{rig-eq: vi to v} and the fact the weak convergence is metrizable as  by Remark \ref{rem: NNNN} we get $\Vert e( (\nabla T^{k}_i)^T \nabla u^k_i) \Vert_{L^2(\Omega)} \le C$ for a constant independent of $k,i$.  Moreover, \eqref{rig-eq: comp2}(ii)  and \eqref{rig-eq: comp3}(i) directly follow from the corresponding estimates for the functions $u^k_i$. Finally, to see \eqref{rig-eq: comp3}(ii) it suffices to prove
\begin{align}\label{new15}
\liminf\nolimits_{i \to \infty} \big({\cal H}^1(J_{u_i} \setminus \partial  P_i) + {\cal H}^1(\partial  P_i\cap \Omega)\big) \ge  {\cal H}^1(J_{u} \setminus \partial  P) + {\cal H}^1(\partial  P\cap \Omega),
\end{align}
where for shorthand $\partial  P_i = \bigcup_j \partial^*  P_{i,j}$ and $\partial  P = \bigcup_j \partial^*  P_{j}$. This can be derived arguing as in \eqref{rig-eq: compa1}:  We may consider an infinitesimal perturbation of the form $\tilde{u}_i = u_i + \sum_j b_j \chi_{P_{i,j}}$,  $\tilde{u} = u + \sum_j b_j \chi_{P_{j}}$ with $b_j$ small such that ${\cal H}^1(\partial P_i \setminus J_{\tilde{u}_i}) ={\cal H}^1(\partial P \setminus J_{\tilde{u}})= 0$  and the convergence in \eqref{rig-eq: vi to v} still holds after replacing $u_i$, $u$ by $\tilde{u}_i$, $\tilde{u}$, respectively.  Then the claim follows from \eqref{rig-eq: vi to v}. Consequently, ${\cal P} \in {\cal Z}_P$ due to Definition \ref{rig-def: ad,coar}(i).

\smallskip

\BBB \emph{Step III:} \EEE It remains to show \eqref{rig-eq: vi to v}. Clearly, we have $\Vert e(\nabla T^T \nabla u_i)\Vert^2_{L^2(\Omega)} \le C$ and ${\cal H}^1(J_{u_i}) \le C$ for a constant independent of $i \in \N$. This follows by a lower semicontinuity argument using \eqref{rig-eq: comp3}(ii) and $\Vert e((\nabla T^{k}_i)^T \nabla u^k_i) \Vert_{L^2(\Omega)} \le C$  by Remark \ref{rem: NNNN}. Thus, in order to apply Theorem \ref{rig-th: GSBD comp}, we have to \BBB select $u_i \in {\cal Z}_u({\cal P}_i)$ suitably and to \EEE find an increasing continuous function $\psi: [0,\infty) \to [0,\infty)$ with $\lim_{\BBB t \EEE \to \infty} \psi(\BBB t \EEE ) = + \infty$ such that $\int_\Omega \psi(|u_i|) \le C$ \BBB independently of $i \in \N$. \EEE

We proceed similarly as in the proof of Theorem \ref{rig-th: comp1} and define $u_i$ iteratively. Choose $u_1 \in {\cal Z}_u({\cal P}_1)$ arbitrarily.   Given $u_i$  we define $u_{i+1}$ as follows.   Consider some $P_{i+1,j}$ and \BBB recalling \eqref{new10} \EEE choose $l_{1,j} <l_{2,j} < \ldots$ such that $P_{i+1,j} = \bigcup^\infty_{k=1} P_{i,l_{k,j}}$ up to an ${\cal L}^2$- negligible set (observe that the union may also be finite). Choose some $\tilde{u}_{i+1} \in {\cal Z}_u({\cal P}_{i+1})$. By Lemma \ref{rig-lemma: comp2.4} for ${\cal P} = {\cal P}_i, \hat{\cal P} = {\cal P}_{i+1}$ we get $(u_i- \tilde{u}_{i+1})\chi_{P_{i+1,j}} = \sum^\infty_{k=1} \BBB \nabla T \EEE(A_{l_{k,j}}\cdot + b_{l_{k,j}})\chi_{P_{i,l_{k,j}}}$ for $A_{l_{k,j}} \in \R^{2 \times 2}_{\rm skew}$, $b_{l_{k,j}} \in \R^2$. Now define 
$$u_{i+1}(x) = \tilde{u}_{i+1}(x) + \BBB \nabla T(x) \EEE (A_{l_{1,j}}\, x + b_{l_{1,j}})$$for $x \in P_{i+1,j}$ and observe that $u_i = u_{i+1}$ on $P_{i,l_{1,j}}$. Proceeding in this way on all $P_{i+1,j}$ we find some $\tilde{A}^{i+1} \in {\cal A}({\cal P}_{i+1})$ such that $u_{i+1}:= \tilde{u}_{i+1} +  \nabla T \tilde{A}^{i+1} \in {\cal Z}_u({\cal P}_{i+1})$ applying Lemma \ref{rig-lemma: comp2.4} for ${\cal P} = \hat{\cal P} = {\cal P}_{i+1}$. Moreover, there is a corresponding $A^i \in {\cal A}({\cal P}_{i})$ such that $u_{i+1} = u_i +  \nabla T  A^i$ with $A^i =0$ on $\bigcup_j P_{i,l_{1,j}}$.

We now show that $\sum_{i \in \N} |A^i(x)| < + \infty$ \BBB for a.e. $x\in\Omega$. \EEE To see this, we recall that \BBB $|P_{i,j} \triangle P_{j}| \to 0$ \EEE for all $j \in \N$. Consequently, as due to the total order of the partitions  the sets $P_{i,j}$ are increasing  for fixed $j \in \N$, the construction of the functions $(u_i)_i$ implies $A^i = 0$ on $P_{i,j}$ for $i$ so large that $ |P_{i,j}| > \frac{1}{2} |P_j|$. Thus, for a.e. $x \in P_{j}$ the sum $\sum_{i \ge 1}  |A^i(x)|$ is a finite  sum and therefore finite. As  $j \in \N$ \BBB was arbitrary, \EEE we obtain  $\sum_{i \in \N} |A^i| < + \infty$ almost everywhere. 

 Therefore, the function $v:= |u_1| + \sum_{l\in \N} |A^l|$ is finite a.e. in $\Omega$ and we apply Lemma \ref{rig-lemma: concave function2} on the sequence $v_k = v\chi_{\lbrace v \le k\rbrace}$ to find an increasing continuous function ${\psi}: [0,\infty) \to [0,\infty)$ with $\lim_{t \to \infty}{\psi}(t) = \infty$ such that by Fatou's lemma $\Vert {\psi}(v) \Vert_{L^1(\Omega)} \le \liminf_{k \to \infty} \Vert {\psi}(v_k) \Vert_{L^1(\Omega)} \le C < \infty$. Using the definition $u_{i+1} = u_i +  \nabla T A^i$ and the monotonicity of ${\psi}$ we find  $\Vert {\psi}(|u_i|) \Vert_{L^1(\Omega)} \le \Vert {\psi}(|u_1| + \sum_{l \in \N} |A^l|) \Vert_{L^1(\Omega)} \BBB \le C \EEE < \infty$ for all $i\in\N$, as desired.  \eop

After these preparatory lemmas we are in a position to prove Theorem \ref{rig-th: comp2}.

\noindent {\em Proof of Theorem \ref{rig-th: comp2}.} First, (i) follows from Lemma \ref{rig-lemma: comp2.5}. The uniqueness of the coarsest partition is a consequence of Lemma \ref{rig-lemma: comp2.2} and Lemma \ref{rig-lemma: comp2.1}. We obtain existence by Zorn's lemma: As $({\cal Z}_P, \ge)$ is a partial order and every chain has an upper bound by Lemma \ref{rig-lemma: comp2.3}, there exists a maximal element $\bar{ \cal P} \in {\cal Z}_P$. Lemma \ref{rig-lemma: comp2.1} shows that $\bar{ \cal P}$ is a coarsest partition which gives (ii). Finally, assertion (iii), namely ${\cal Z}_u(\bar{\cal P}) = v + \nabla T {\cal A}(\bar{\cal P})$ for some $v \in {\cal Z}_u(\bar{\cal P})$, follows from Lemma \ref{rig-lemma: comp2.4} for the choice ${\cal P} = \hat{\cal P} = \bar{\cal P}$. \eop

%----------------------------------------------------------------------------------------------------------------------
\section{The effective linearized Griffith model}\label{rig-sec: gamma}
%-----------------------------------------------------------------------------------------------------------------------

In this final section we identify the effective linearized Griffith functional via $\Gamma$-convergence and derive a cleavage law for the limiting model.  

%----------------------------------------------------------------------------------------------------------------------
\subsection{Derivation of linearized models via  $\Gamma$-convergence}\label{rig-sec: sub, gamma1}
%-----------------------------------------------------------------------------------------------------------------------

We now give the proof of Theorem \ref{rig-th: gammaconv}. 

\noindent {\em Proof of Theorem \ref{rig-th: gammaconv}.} (i) Thanks to the preparations in the last section the lower bound is almost immediate. Let $(u,{\cal P},T) \in {\cal D}_\infty$ be given as well as a sequence $(y_k)_k \subset SBV_M(\Omega)$ with \BBB $y_k \to (u,{\cal P},T)$, i.e. by Definition \ref{def:conv} the are triples \EEE $(u_k,{\cal P}_k,T_k) \in {\cal D}$ such that \eqref{rig-eq: comp2}-\eqref{rig-eq: comp3} hold. Due to \eqref{rig-eq: comp3}(ii) it suffices to show
$$\liminf_{k \to \infty} \frac{1}{\eps_k}\int_\Omega W(\nabla y_k) \ge \int_\Omega \frac{1}{2} Q(e(\nabla T^T \nabla u) ).$$
We proceed as in \eqref{rig-eq: ff}:  Recall that $W(G) = \frac{1}{2}Q( e(G-\Id)) + \omega_{\BBB W \EEE } (G-\Id)$ with $\sup\lbrace |F|^{-3}\omega_W(F): |F| \le  1\rbrace \le C$ by the assumptions on $W$, where $Q  =  D^2 W(\Id)$. We compute  by \eqref{rig-eq: comp3}(i)
\begin{align*}
\frac{1}{\eps_k} \int_{\Omega} W(\nabla y_k ) &  \ge \frac{1}{\eps_k} \int_{\Omega}  W(\Id + \sqrt{\eps_k}\nabla T_k^T \nabla u_k ) + o(1)\\
& =  \int_{\Omega}  \frac{1}{2}\Big( Q(e(\nabla T_k^T \nabla u_k)) + \frac{1}{\eps_k}\omega_W(\sqrt{\eps_k} \nabla T_k^T\nabla u_k)  \Big) +  o(1)
\end{align*}
as $k \to \infty$. The second term converges to $0$ arguing as in \eqref{rig-eq: ff2} and using $\Vert \nabla u_k\Vert_\infty \le c\eps_k^{-1/8}$ (see \eqref {rig-eq: comp2}(ii)). As $e( \nabla T_k^T\nabla u_k) \rightharpoonup e(\nabla T^T \nabla u)$ weakly in $L^2(\Omega,\R^{2\times 2}_{\rm sym})$ by \eqref{rig-eq: comp1-2}(ii) and $Q$ is convex, we conclude 
$$\liminf_{k \to \infty} \frac{1}{\eps_k} \int_{\Omega} W(\nabla y_k)\ge \int_{\Omega} \frac{1}{2} Q(e(\nabla T^T \nabla u)).$$

(ii) By a general density result in the theory of $\Gamma$-convergence together with Theorem \ref{rig-th: cortesani2}  and the fact that the limiting functional $E(u,{\cal P}, T)$ is continuous in $u$ with respect to the convergence given in Theorem \ref{rig-th: cortesani2}, it suffices to provide recovery sequences for functions $u$ with $u \in W^{1,\infty}(\Omega \setminus \overline{J_u})$, where $J_u$ is contained in the union of a finite number of closed connected pieces of $C^1$- curves. Moreover,  as in the proof of Theorem \ref{rig-th: comp1} \BBB (see paragraph before \eqref{rig-eq: compa1}) \EEE we may assume that ${\cal H}^1( (\bigcup_j \partial^* P_j \cap \Omega) \setminus J_{u}  ) = 0$ up to an infinitesimal small perturbation  of $u$ (a similar argument was used below \BBB \eqref{new15}). \EEE  

Let $(u,{\cal P},T) \in {\cal D}_\infty$ and $\eps_k \to 0$ be given. Define $y_k(x) = T(x) + \sqrt{\eps_k} u(x)$ for all $x \in \Omega$. It is not hard to see that $(y_k)_k \subset SBV_M(\Omega)$ for  $\eps_k$ small enough \BBB (and $M$ not too small). \EEE Moreover, define ${\cal P}_k = {\cal P}$, $T_k= T$ and  $u_k =   \eps_k^{-1/2} \big( y_k - T_k\big) \equiv u$ for all $k \in \N$.  Then \eqref{rig-eq: comp2}(i) and \eqref{rig-eq: comp1}-\eqref{rig-eq: comp3}  hold trivially. To see  \eqref{rig-eq: comp2}(ii), it suffices to note that $\Vert \nabla u_k \Vert_\infty = \Vert \nabla u \Vert_\infty \le C \le C\eps_k^{-1/8}$. Consequently, $y_k \to (u, {\cal P},T)$ in the sense of Definition  \ref{def:conv}. \EEE

We finally confirm $\lim_{k \to \infty} E_{\eps_k}(y_k) = E(u,{\cal P},T)$. As clearly   $\lim_{k \to \infty}{\cal H}^1(\BBB J_{y_k} \EEE ) =  {\cal H}^1(\bigcup_j\partial^* P_j\cap \Omega) + {\cal H}^1(J_u \setminus \bigcup_j\partial^* P_j)$,  it suffices to show $\lim_{k \to \infty} \frac{1}{\eps_k}\int_\Omega W(\nabla y_k) = \int_\Omega \frac{1}{2} Q(e(\nabla T^T \nabla u))$. Using again that $W(G) =  \frac{1}{2}Q( e(G-\Id)) + \omega_W(G-\Id)$ \BBB and the frame indifference of $W$ \EEE we compute
\begin{align}\label{new22}
\begin{split}
\frac{1}{\eps_k} \int_{\Omega} W(\nabla y_k ) & =\frac{1}{\eps_k} \int_{\Omega} W( \BBB \Id + \sqrt{\eps_k} \nabla T_k^T \nabla u_k \EEE )\\
& =  \int_{\Omega}  \Big( \frac{1}{2} Q(e(\nabla T_k^T \nabla u_k)) + \frac{1}{\eps_k}\omega_W(\sqrt{\eps_k} \nabla T_k^T\nabla u_k)  \Big)  \\ &= \int_{\Omega}   \frac{1}{2} Q(e(\nabla T^T \nabla u)) + O(\sqrt{\eps_k})\to \int_{\Omega}   \frac{1}{2} Q(e(\nabla T^T \nabla u)).
\end{split}
\end{align}
This concludes the proof. \eop

The proof of Corollary \ref{rig-cor: gamma} is now straightforward. 

\noindent {\em Proof of Corollary \ref{rig-cor: gamma}.} To see the liminf-inequality, assume  $y_{\eps_k} \to y$ in $L^1(\Omega)$ for $k \to \infty$ and without restriction that   $E_{\eps_k}(y_{\eps_k}) \le C$. \BBB By Theorem \ref{rig-th: comp1}   we find a limiting triple $(u,{\cal P},T) \in {\cal D}_\infty$ such that $y_{\eps_k} \to (u,{\cal P},T)$ in the sense of Definition \ref{def:conv}. \EEE By \eqref{rig-eq: comp2}(i), \eqref{rig-eq: comp1}(ii) we obtain $y = T$. \BBB As $T \in {\cal R}({\cal P})$, we get $\mathcal{H}^1(J_T) \le \mathcal{H}^1(\bigcup_j \partial^* P_j \cap \Omega)$, where ${\cal P} = (P_j)_j$. \EEE Thus, Theorem \ref{rig-th: gammaconv}(i) yields $\liminf_{k \to \infty} E_{\eps_k}(y_k) \BBB \ge E(u,{\cal P},T) \ge {\cal H}^1(J_T) \EEE = E_{\rm seg}(y)$. A recovery sequence is obviously given by $y_k = y$ for all $k \in \N$. \eop

We close this section with the proof of Lemma \ref{lemmanew} and Theorem \ref{rig-th: gammaconv2}.

\BBB  

\noindent {\em Proof of Lemma \ref{lemmanew}.} Consider triples $(u_k,{\cal P}_k,T_k)$ and  $(g_k,{\cal P}^g_k,T^g_k)$   such that the triples converge to $(u, {\cal P},T)$ and $(g,{\cal P}_g,T_g)$, respectively, in the sense of \eqref{rig-eq: comp2}-\eqref{rig-eq: comp3}. Since $\eps_k^{-1/2}(y_k - f_k)$ is bounded in $L^2(\Omega)$ by \eqref{ennew}, $\eps_k^{-1/2}(T_k - T_k^g) + (u_k - g_k) - \eps_k^{-1/2}(y_k - f_k)$ converges  a.e.  by \eqref{rig-eq: comp2}(i), and $u_k - g_k$ converges a.e.  by \eqref{rig-eq: comp1-2}(i), we get $\eps_k^{-1/2}(T_k - T_k^g)$ converges (up to a subsequence) in measure on $\Omega$. This implies $T = T_g$ by \eqref{rig-eq: comp1}(ii). 

Moreover, suppose ${\cal P}_g \ge {\cal P}$ was wrong (recall \eqref{new10}). We may then assume after reordering that $P_1^g \cap {P}_1$ and $P_2^g \cap {P}_1$ have positive ${\cal L}^2$-measure. Since $\eps_k^{-1/2}(T_k - T_k^g)$ converges in measure on $P_1$, we can  argue exactly as in the proof of  Lemma \ref{rig-lemma: comp2.1}(2)  to see that the partition ${\cal P}_g$ is not coarsest, which contradicts the assumption.  \eop

\noindent {\em Proof of Theorem \ref{rig-th: gammaconv2}.} First consider the lower bound.  Let $(u,{\cal P},T) \in {\cal D}_\infty$ be given as well as a sequence $(y_k)_k \subset SBV_M(\Omega)$ with $y_k \to (u,{\cal P},T)$. If $(u,{\cal P},T) \notin {\cal C}_g$,  Lemma \ref{lemmanew} implies $\liminf_{k \to \infty} F_{\eps_k}(y_k) = \infty$. If $(u,{\cal P},T) \in {\cal C}_g$, we recall Theorem \ref{rig-th: gammaconv}(i) and see that it suffices to show 
\begin{align}\label{new21}
\liminf_{k \to \infty}\eps_k^{-1}  \Vert y_k - f_k \Vert^2_{L^2(\Omega)} \ge  \min_{v \in u + \nabla T \mathcal{A}(\mathcal{P})}  \Vert v- g \Vert^2_{L^2(\Omega)}.
\end{align}
Consider the  triples $(u_k,{\cal P}_k,T_k)$ and  $(g_k,{\cal P}^g_k,T^g_k)$  as in the previous proof. Recall $\eps_k^{-1/2}(y_k - f_k ) - \eps_k^{-1/2}(T_k - T_k^g) - (u_k - g_k) \to 0$ a.e. by \eqref{rig-eq: comp2}(i) and $u_k - g_k \to u - g$ a.e. by \eqref{rig-eq: comp1-2}(i). Now it is enough to show that for each $P_j$ there are $A_j \in \R^{2 \times 2}_{\rm skew}$
 and $d_j \in \R^2$ such that $\eps_k^{-1/2}(T_k - T_k^g) \to \nabla T(A_j\cdot+d_j)$ a.e. on $P_j$. Then \eqref{new21} follows from  Fatou's lemma. We recall from the proof of  Lemma \ref{lemmanew} that $\eps_k^{-1/2}(T_k - T_k^g)$ converges in measure on $P_j$ and that $P_j  \subset P_i^g$ for some $P_i^g$. In particular, this implies $|R^k_j - R_i^{g,k}| +  |b^k_j - b_i^{g,k}| \le C\sqrt{\eps_k}$, where $R_i^{g,k}\cdot + b_i^{g,k}$ denotes the rigid motion associated to $T_k^g$. Repeating the argument in \eqref{rig-eq: Tdef} we obtain the desired convergence. 
 
For the construction of recovery sequences we mainly follow  Theorem \ref{rig-th: gammaconv}(ii) and only   indicate the necessary adaptions. Let $(u,{\cal P},T) \in {\cal C}_g$ be given with $u$ having the specific regularity assumed in Theorem \ref{rig-th: gammaconv}(ii), particularly $\nabla u \in L^\infty(\Omega)$. Let $A \in {\cal A}({\cal P})$ such that $v:= u + \nabla T A$ realizes the minimum in \eqref{new21}. 

As $u,g \in L^2(\Omega)$ and thus $\nabla T A \in L^2(\Omega)$, we can choose  a sequence $(A^k)_k \subset {\cal A}({\cal P}_k^g) $ such that $\nabla T_k^g A^k \to \nabla T_g A = \nabla T A$ in $L^2(\Omega)$ and $\sqrt{\eps_k}|\nabla A^k|^2   \to 0$ in $L^2(\Omega)$. Select $T_k =  \sum_j (R^{k}_j\cdot+b^{k}_j)\chi_{P^{g,k}_j} \in {\cal R}({\cal P}_k^g)$ such that (cf. before \eqref{new5} for a similar construction)
$$\dist(\nabla (T_k^g + \sqrt{\eps_k}\nabla T_k^g A^k) ,SO(2)) = |\nabla (T_k^g + \sqrt{\eps_k}\nabla T_k^g  A^k) - \nabla T_k| \text{ on } \Omega$$
and  on each component the translations of $T_k$ and $T_k^g + \sqrt{\eps_k}\nabla T_k^g A^k$ coincide, i.e. 
$$b^{k}_j = T_k(x) - \nabla T_k(x)\,x =  (T_k^g + \sqrt{\eps_k}\nabla T_k^g A^k)(x) - \nabla (T_k^g + \sqrt{\eps_k}\nabla T_k^g A^k)(x) \, x$$
 for all $x \in P^{g,k}_j$ and all $j \in \N$. Note that using \eqref{rig-eq: linearization} a short calculation implies 
\begin{align}\label{new25}
| T_k^g + \sqrt{\eps_k}\nabla T_k^g A^k - T_k|  + |\nabla T_k^g + \sqrt{\eps_k}\nabla T_k^g \nabla A^k - \nabla T_k|  \le C\eps_k|\nabla A^k|^2 
\end{align}
pointwise a.e. in $\Omega$. We define $y_k = T_k +  \sqrt{\eps_k} u$,   $u_k = u$  and ${\cal P}_k = {\cal P}$ for all $k \in \N$. We note that $T_k \in {\cal R}({\cal P}_k)$ since $T_k \in {\cal R}({\cal P}^g_k)$,  ${\cal P}_k \le  {\cal P}_g$ by $(u,{\cal P},T) \in {\cal C}_g$ and $ {\cal P}_k^g = {\cal P}^g$ by assumption. Moreover,  $\nabla T_k^T \nabla y_k = \Id + \sqrt{\eps_k}\nabla T_k^T \nabla u$ and by \eqref{new25}
\begin{align}\label{new23}
\nabla T_k^T \nabla u = (\nabla T_k^g)^T \nabla u + O(\Vert \nabla u\Vert_\infty \,\eps_k|\nabla A^k|^2) + O(\Vert \nabla u\Vert_\infty \sqrt{\eps_k}|\nabla A^k|).
\end{align}  
Using $\sqrt{\eps_k}|\nabla A^k|^2  \to 0$ in $L^2(\Omega)$ and $T_k^g \to T, \nabla T_k^g \to \nabla T$ in $L^2(\Omega)$, one can check that  $(y_k)_k \subset SBV_M(\Omega)$ and that \eqref{rig-eq: comp2}-\eqref{rig-eq: comp3} hold, i.e. $y_k \to (u,{\cal P},T)$ (see also proof of Theorem \ref{rig-th: gammaconv}). Likewise, \eqref{new23} together with the calculation in \eqref{new22} shows that $\lim_{k\to \infty} E_{\eps_k}(y_k) = E(u,{\cal P},T)$. Finally, recalling $\eps_k^{-1/2}(f_k - T_k^g) \to g$ in $L^2(\Omega)$ (see before Lemma \ref{lemmanew}), $\nabla T_k^g A^k \to \nabla T A$ in $L^2(\Omega)$ and \eqref{new25}  we conclude  
\begin{align*}
\eps_k^{-1/2}\Vert y_k - f_k \Vert_{L^2(\Omega)} &=  \Vert \eps_k^{-1/2}  (T_k - T_k^g) +   u - \eps_k^{-1/2}(f_k - T_k^g)\Vert_{L^2(\Omega)}\\& =
\Vert \nabla T_k^g A^k +    u - \eps_k^{-1/2}(f_k - T_k^g) \Vert_{L^2(\Omega)} +  O(\Vert \eps_k^{1/2}|\nabla A^k|^2 \Vert_{L^2(\Omega)})\\&\to \Vert \nabla T A +    u - g\Vert_{L^2(\Omega)} =\Vert v - g\Vert_{L^2(\Omega)}.
\end{align*}
Finally, the convergence result for minimum problems and minimizers follows from a general result in the theory of $\Gamma$-convergence (see \cite[Chapter 7]{DalMaso:93}). \eop
 \EEE

%----------------------------------------------------------------------------------------------------------------------
\subsection{An application to cleavage laws}\label{rig-sec: sub, gamma2}
%-----------------------------------------------------------------------------------------------------------------------

We are finally in a position to prove the cleavage law in Theorem \ref{rig-th: cleavage-cont}. Analogous results for the case of expansive boundary values have been obtained in \cite{Mora:2010} and \cite{FriedrichSchmidt:2014.2}. We thus do not repeat all the steps of these proofs but rather concentrate on the additional arguments necessary in our general setting (see \eqref{rig-eq: Griffith en-lim}) in which we particularly can extend the aforementioned results to the case of compression. 
\smallskip

\noindent {\em Proof of Theorem \ref{rig-th: cleavage-cont}.} Let $(y_{\eps_k})_{k}$ be a sequence of almost minimizers. Passing to a suitable subsequence, by Theorem \ref{rig-th: comp1} we obtain a triple $(u_k, {\cal P}_k, T_k) \in {\cal D}$ and a limiting triple $(u, {\cal P}, T) \in {\cal D}_\infty$ such that \eqref{rig-eq: comp2}-\eqref{rig-eq: comp3} hold and
$$ E(u, {\cal P}, T) \le \liminf\nolimits_{\eps \to 0} \inf  \lbrace E_\eps(y): y \in {\cal A}(a_\eps)\rbrace$$
by Theorem \ref{rig-th: gammaconv}(i). \BBB Write $T_k = \sum_j (R_j^k\cdot + b_j^k)\chi_{P_j^k}$ and ${\cal P} =(P_j)_j$. \EEE Due to the boundary conditions \BBB and \eqref{rig-eq: comp2}(i),\eqref{rig-eq: comp1-2}(i), \EEE  on each component $P_j \in {\cal P}$ we find $A_j \in \R^{2 \times 2}_{\rm skew}$ and $b_j \in \R^2$ such that 
\begin{align}\label{rig-eq: bc-cleav2}
\begin{split}
u_1(x) &= \lim\nolimits_{k \to \infty} \eps_k^{-1/2}(  \e_1 \cdot (\Id - R^k_j)\, x  - \e_1 \cdot b^k_j+  a_{\eps_k} x_1 ) \\ &= \e_1 \cdot A_j\, x  +\e_1 \cdot b_j +  a x_1
\end{split}
\end{align}
for a.e. $x \in \Omega'$  with $x_1 <0$ or $x_1 > l$ and $x \in P_j$. In particular, this implies
\begin{align}\label{rig-eq: bc-cleav}
 u_1(x_1,x_2) - u_1(\hat{x}_1,x_2)   = |x_1 - \hat{x}_1| a 
\end{align}
for a.e. $x \in  \Omega'$  with  $\hat{x}_1 < 0$, $x_1 > l$ and $(x_1,x_2), (\hat{x}_1,x_2) \in P_j$.

We first derive the limiting minimal energy and postpone the characterization of the sequence of almost minimizers to the end of the proof. The argument in \eqref{rig-eq: bc-cleav2} shows that $\nabla T = \Id$ on $P_j$ if $|P_j \cap \lbrace x: x_1 < 0 \text{ or } x_1 > l \rbrace|>0$.  It is not restrictive to  assume $\nabla T^T \nabla u = \nabla u$ a.e.  Indeed, we may replace $u$ by \BBB $0$ \EEE in a component $P_j$ which does not intersect the boundaries  without \BBB increasing \EEE the energy. By \eqref{rig-eq: Griffith en-lim}, a slicing argument in $GSBD$ (see Theorem \ref{clea-th: slic}) and the fact that $\inf \lbrace Q(F): \e_1^T F \e_1 = a\rbrace = \alpha a^2$ (see Section \ref{rig-sec: sub, cleavage}) we obtain
\begin{align}\label{eq:E-uS}
\begin{split}
 E(u, {\cal P}, T) &\ge \int_{\BBB \Omega \EEE }  \frac{1}{2} Q(e(\nabla u)) + \int_{J_u} | \BBB \nu_u \EEE \cdot \e_1| d{\cal H}^1 + {\cal E}(u) \\
&\ge  \int_0^1 \Big( \int_0^l \frac{\alpha}{2} (\e_1^T \nabla u(x) \e_1)^2 \, dx_1 + S^{x_2}(u) \Big) \, dx_2 + {\cal E}(u),
\end{split}
\end{align} 
where $S^{x_2}$ denotes the number of jumps of $u_1$ on a slice $(-\eta,l+\eta) \times \lbrace x_2 \rbrace$ and ${\cal E}(u) = \int_{J_u} (1 -|\nu_u \cdot \e_1|) d{\cal H}^1$. If $S^{x_2} \ge 1$, the inner integral is bounded from below by $1$. By the structure theorem for Caccioppoli partitions (see Theorem \ref{th: local structure}) we find that $((-\eta ,0) \cup (l, l +\eta)) \times \lbrace x_2 \rbrace \subset P_j$   for some $j \in \N$ for ${\cal H}^1$-a.e. $x_2$ with $S^{x_2} = 0$. Consequently, if $\#S^{x_2} =0$, by applying Jensen's inequality we derive that the term is bounded from below by $\frac{1}{2} \alpha l a^2$ due to the boundary conditions \eqref{rig-eq: bc-cleav}. This implies $E(u) \ge \min \lbrace \frac{1}{2} \alpha l a^2,1 \rbrace$. 

Otherwise, it is not hard to see that the configurations $y^{\rm el}_{\eps_k}(x) = x + F^{a_{\eps_k}} \, x$ for $x \in \Omega'$  satisfy $E_{\eps_k}(y^{\rm el}_{\eps_k}) \to \frac{1}{2}\alpha l a^2$ for $k\to \infty$. Likewise, we get $E_{\eps_k}(y^{\rm cr}_{\eps_k}) = 1$ for all $k \in \N$, where   $y^{\rm cr}_{\eps_k}(x) = x \chi_{x_1 < \frac{1}{2}} + (x + (l a_{\eps_k},0)) \chi_{x_1 > \frac{1}{2}}$ for $x \in \Omega$ and $y^{\rm cr}_{\eps_k} = (x_1(1+a_{\eps_k}) , x_2 )$ for $x \in \Omega' \setminus \Omega$. This completes \eqref{rig-eq: cleavage en}.

\smallskip

It remains to characterize the sequences of almost minimizers. \BBB Let $(y_{\eps_k})_{k}$ be a sequence of almost minimizers 
and $(u,{\cal P},T) \in {\cal D}_\infty$ a limiting triple as considered before \eqref{rig-eq: bc-cleav2}. Again we may suppose $\nabla T^T \nabla u = \nabla u$ a.e. \EEE We let first $|a| < a_{\rm crit}$ and  follow the arguments in the proof of \cite[Theorem 2.4]{FriedrichSchmidt:2014.2}. Since $E(u,{\cal P},T) = \frac{1}{2}\alpha la^2$,  we infer from \eqref{eq:E-uS} that $u$
has no jump on a.e.\ slice $(-\eta,l+\eta) \times \left\{x_2\right\}$ and satisfies $\e_1^{T} \nabla u \, \e_1 = a$ a.e.\ by the imposed boundary values and the strict convexity of the mapping $t \mapsto t^2$ on \BBB $\R$. \EEE Thus, if $J_{u} \neq
\emptyset$,  a crack normal must satisfy $\nu_{u} = \pm \e_2$ ${\cal H}^1$-a.e.
Taking additionally ${\cal E}(u)$ into account we find $J_{u} = \emptyset$ up to an ${\cal H}^1$ negligible set, i.e., $u \in H^1(\Omega')$.  By the strict convexity of $Q$ on symmetric
matrices and the boundary values \eqref{rig-eq: bc-cleav2} we see that the derivative has the form
$$ \nabla u(x) = \begin{footnotesize} F^a + A \end{footnotesize} \text{ for a.e.\ $x \in \Omega$} $$
for a suitable $A \in \R^{2 \times 2}_{\rm skew}$. Since $\Omega$ is connected, we conclude 
$$u(x) = F^a \,x + A\,x + b$$ 
for $x \in \Omega$  and some  $b \in \R^2$.  In particular, this implies ${\cal P}$ consists only of $P_1 = \Omega'$ and thus by \eqref{rig-eq: bc-cleav2} we get $A = \lim_{k \to \infty} \eps_k^{-1/2}(\Id - R^k_1)$ and $\e_1 \cdot b = - \lim_{k \to \infty} \eps^{-1/2}_{k}\e_1 \cdot b^k_1$. Let $ s=   \lim_{k \to \infty} \e_2 \cdot (\eps_k^{-1/2} b^k_1 +b)$, which exists by \eqref{rig-eq: comp2}(i), \eqref{rig-eq: comp1-2}(i). We now conclude by \eqref{rig-eq: comp2}(i), \BBB \eqref{rig-eq: comp1-2}(i) \EEE for  a.e. $x \in \Omega$ 
\begin{align}\label{rig-eq: lincleav}
\begin{split}
\BBB \bar{u}(x):= \EEE \lim_{k\to \infty} \eps_k^{-1/2}(y_{\eps_k}(x) - x)  & = u(x) + \lim_{k\to \infty} \eps_k^{-1/2} \big((R^k_1 - \Id)\,x +  b^k_1 \big) \\
& = u(x) - A\,x - b + (0,s)= (0,s) + F^a \, x,
\end{split}
\end{align}
\BBB i.e. $\bar{u}$ fulfills Theorem \ref{rig-th: cleavage-cont}(i). \EEE If $|a| > a_{\rm crit}$, we again consider the lower bound  \eqref{eq:E-uS} and now obtain that on
a.e.\ slice $(0,l) \times \left\{x_2\right\}$ a minimizer $u$ has  precisely one jump
and that $\e_1^{T} \nabla u \, \e_1 = 0$ a.e. By the strict convexity of $Q$ on symmetric matrices we then derive that 
$\nabla u$ is antisymmetric a.e.  As a consequence, the linearized piecewise rigidity estimate for $SBD$ functions  (see \cite[Theorem A.1]{Chambolle-Giacomini-Ponsiglione:2007} or the remark below Theorem \ref{rig-cor: cgp}) yields that there is a Caccioppoli partition $(E_i)_i$ of $\Omega$ such that 
$$ u(x) = \sum\nolimits_{i} (A_i x + b_i) \chi_{E_i}(x) \quad \text{and} \quad 
   J_u = \bigcup\nolimits_{i} \partial^* E_i \cap \Omega, $$
where $A_i\in \R^{2 \times 2}_{\rm skew}$ and $b_i \in \R^2$. (Note that indeed the linearized rigidity estimate can also be applied in the $GSBD$-setting as it relies on a slicing argument and an approximation which is also available in the generalized framework, see \cite[Section 3.3]{Iurlano:13}. The only difference is that the approximation does not converge in $L^1$ but only pointwise a.e. which does not affect the argument.) 

As ${\cal E}(u) = 0$, we also note that $\nu_u = \pm \e_1$ a.e.\ on $J_u$. Following the arguments in \cite{Mora:2010}, in particular using regularity results for boundary curves of sets of finite perimeter and exhausting the sets $\partial^* E_i$ with Jordan curves, we find that 
$$ J_u = \bigcup\nolimits_{i} \partial^* E_i  \cap \Omega \subset (p, 0) + \R \e_1 $$ 
for some $p \in (0,l)$. We thus obtain that $(E_i)_i$ consists of only two sets and $u$ has the form
$$u(x) = \begin{cases} A_1 \, x + b_1 & \text{for } x_1 < p, \\ A_2 \, x + b_2 & \text{for } x_1 > p, \end{cases} $$
for $A_i \in \R^{2 \times 2}_{\rm skew}$and $b_i \in \R^2$, $i=1,2$.  Now repeating the calculation in \eqref{rig-eq: lincleav} for the sets $P_1 = \lbrace x \in \Omega': x_1 < p\rbrace$ and $P_2 = \Omega' \setminus P_1$ we find  $s, t \in \R$ such that for $x \in \Omega$ a.e. 
 \begin{align*}
\BBB \bar{u}(x):= \EEE \lim_{k\to \infty} \eps_k^{-\frac{1}{2}}(y_{\eps_k}(x) - x)  & = u(x) - (A_1\,x + b_1)\chi_{x_1 < p}(x) - (A_2\,x + b_2)\chi_{x_1> p}(x) \\ & \ \ \ + (0,s)\chi_{x_1 < p}(x) + ((l a,t))\chi_{x_1 > p}(x).
\end{align*}
\BBB Then $\bar{u}$ satisfies Theorem \ref{rig-th: cleavage-cont}(ii). \EEE This concludes the proof. \eop

\textbf{Acknowledgements}  I am very grateful to Bernd Schmidt for many stimulating discussions and valuable comments from which the results of this paper and their exposition have benefited a lot. \BBB Moreover, I am gratefully indebted to the referee for her/his careful reading of the manuscript and many helpful suggestions. \EEE

%--------------------------------------------------------------------------

%--------------------------------------------------------------------------
%--------------------------------------------------------------------------
 \typeout{References}

\end{document}